\documentclass[a4paper,12pt,draft]{amsart}

\usepackage[latin1]{inputenc}
\usepackage{bbm}
\usepackage{amssymb}  
\usepackage{amsthm}  

\allowdisplaybreaks

\theoremstyle{plain}
\newtheorem{thm}{Theorem}[section]
\newtheorem{lemma}[thm]{Lemma}
\newtheorem{cor}[thm]{Corollary}

\theoremstyle{definition}
\newtheorem{defn}[thm]{Definition}

\theoremstyle{remark}

\newtheorem{remark}[thm]{Remark}
\newtheorem{example}[thm]{Example}

\newenvironment{enum_i}{%
  \begin{enumerate}%
  }%
  {\end{enumerate}%
}

  {\end{enumerate}%
}

\numberwithin{equation}{section}

\newcommand{\R}{\mathbbm{R}}
\newcommand{\Z}{\mathbbm{Z}}

\newcommand{\intR}{\int_{\mathbbm{R}}}

\newcommand{\sumZ}[1][k]{\displaystyle\sum_{#1 \in \Z}}
\newcommand{\norm}[2]{\left|\!\left|#1\right|\!\right|_{#2}}
\newcommand{\half}{\frac{1}{2}}

\newcommand{\G}{\mathcal{G}}
\newcommand{\Pcal}{\mathcal{P}}

\newcommand{\sgn}{\operatorname{sgn}}

\newcommand{\E}{\mathbbm{E}}
\newcommand{\mesh}[1][]{\left|\pi_{#1}\right|}

\newcommand{\Prob}{\mathbbm{P}}

\newcommand{\comment}[1]{}

\begin{document}

\title[Convergence Rates for SDEs]{Convergence Rates for Approximations of Functionals of SDEs  
}

\author[R. Avikainen]{Rainer Avikainen}
\address{Department of Mathematics and Statistics, P.O. Box 35 (MaD), FIN-40014 University of Jyv\"askyl\"a, Finland} 
\email{avikain@maths.jyu.fi}

\begin{abstract}
We consider upper bounds for the approximation error $\E|g(X)-g(\hat X)|^p$, where $X$ and $\hat X$ are random variables such that $\hat X$ is an approximation of $X$ in the $L_p$-norm, and the function $g$ belongs to certain function classes, which contain e.g. functions of bounded variation. We apply the results to the approximations of a solution of a stochastic differential equation at time $T$ by the Euler and Milstein schemes. For the Euler scheme we provide also a lower bound.
\end{abstract}

\thanks{The author was supported in part by the Finnish Graduate School in Stochastics.}

\maketitle

\section{Introduction}
\label{SectionIntroduction}

Convergence schemes for the solutions of SDEs are rather well known. Let $X$ be the solution of the one-dimensional equation
$$dX_t =  \sigma(t,X_t)\,dW_t + b(t,X_t)\,dt, \ X_0 = x_0,$$
where $W$ is a standard Brownian motion, $t \in [0,T]$, and $\sigma$ and $b$ satisfy certain assumptions specified in Section \ref{SectionAssumptions}.
P.E. Kloeden and  E. Platen \cite{KP} have showed that any order of strong convergence can be achieved by the strong It\^o-Taylor approximations, i.e. for any order $\gamma > 0$ there exists a scheme $X^\pi$ corresponding to a partition $\pi$ of the interval $[0,T]$ with mesh size $\mesh$ such that
$$\E\left(\sup_{0 \leq t \leq T} |X_t-X^\pi_t|\right) < C \mesh^\gamma.$$
The simplest and most common examples are the Euler scheme $X^E$ and the Milstein scheme $X^M$, which have the order of strong convergence 0.5 and 1, respectively. 

Lately this topic has been considered by N. Hofmann, T. Müller-Gronbach and K. Ritter in \cite{HMR1, HMR2},  Hofmann and Müller-Gronbach in \cite{HM}, and Müller-Gronbach in \cite{M2,M3}. In these papers they cover errors with respect to both global and pointwise error criterions. The latest result concerning the pointwise error is due to Müller-Gronbach \cite{M3}, where the author defines certain classes of convergence schemes and finds optimal (adaptive) schemes for each class. 

The pointwise weak error $\E f(X_T) - \E f(X^\pi_T)$ is also considered by several authors, e.g. Kloeden and Platen \cite{KP}, under certain smoothness conditions on $f$. In the case of the Euler scheme the smoothness conditions were relaxed by V. Bally and D. Talay in \cite{BT1}, 
where $f$ was required to be only measurable and bounded. This was done under a Hörmander type condition for the infinitesimal generator of $X$. A recent contribution to this field is given by Moon et al. in \cite{MSTZ}.

This paper considers the error of the type $\E|g(X_T)-g(X^\pi_T)|^p$. The question is motivated by discretization schemes for BSDEs. The terminal condition $g(X_T)$ is approximated by $g(X^\pi_T)$, and $L_p$-estimates for the difference $g(X_T)- g(X^\pi_T)$ are required.
 If $g$ is Lipschitz, this error returns immediately to the error of the underlying scheme. Therefore the aim of this paper is to give results for relevant non-Lipschitz functions. It is proved that if approximations $(X_t^\pi)_{t \in [0,T]}$ satisfy 
 $$\norm{X_T-X^\pi_T}{p} \leq C_p^1\mesh^\gamma$$
 for some $\gamma > 0$ and all $1 \leq p < \infty$, then 
$$\E|g(X_T)-g(X^\pi_T)|^p \leq C_p^2\mesh^{\gamma-\varepsilon}$$
 for any $0 < \varepsilon < \gamma$ and for any $g$ in a 
special class of functions. This class contains functions of bounded variation, polynomials and jump functions controlled by the tail distributions of $X_T$ and $X_T^\pi$, and therefore by Minkowski's inequality the result is true for any sum of these three types of functions. 

Finally, for the Euler scheme a lower bound is presented indicating that the error under consideration can not converge faster than $\mesh^{1/2}$. This is done by giving an example, namely the geometric Brownian motion, for which the lower bound is obtained. There is still a difference, although arbitrarily small, between the upper and lower bounds, and it remains open whether the rates could be equalized. These results are achieved under certain conditions on the SDE, including the existence of a bounded density for the solution $X_T$. 

The structure of the paper is such that the assumptions that hold throughout the paper are collected in Section \ref{SectionAssumptions}. Sections \ref{SectionErrorIndic} -- \ref{SectionExtension} contain the main results concerning upper bounds.  The first result is given in Section \ref{SectionErrorIndic}, where Theorem \ref{IndicatorError} clarifies the convergence rate for indicator functions. This is then applied to the Euler and Milstein schemes in Theorem \ref{IndicatorRateThmEulerMilstein}. The result is then extended to functions of bounded variation in Theorem \ref{BVError} in Section \ref{SectionErrorBV}, and applied to the Euler and Milstein schemes in Theorem \ref{BVRateThmEulerMilstein}. Another extension is developed in Section \ref{SectionExtension}, where the result for the function class $\G_{p, \varphi}$ is given in Theorem \ref{Gpthm}. The class $\G_{p. \varphi}$ is then analyzed in Section \ref{SectionGpAnalysis}, including the result that it contains all polynomials in Theorem \ref{polynomialsinG}. 
An application to the Euler and Milstein schemes is presented in Corollary \ref{ExtensionRateCorEulerMilstein}. Section \ref{SectionLowerBound} contains a lower bound for the convergence, stated in Theorem \ref{lowerboundthm}.   Finally, a proof of a Theorem from the book of Bouleau and L\'epingle \cite{BL} is presented with explicit constants in Appendix \ref{SectionAppendix}.

\section{Assumptions}
\label{SectionAssumptions}
We fix a terminal time $T>0$ and suppose that $(W_t)_{t \in [0,T]}$ is a standard one-dimensional Brownian motion defined on a complete filtered probability space $(\Omega,\mathcal{F},\Prob,(\mathcal{F}_t)_{t \in [0,T]})$, where the filtration is the augmentation of the natural filtration of $W$ and $\mathcal{F}=\mathcal{F}_T$.

We consider a diffusion process $X$, which is a solution to
\begin{equation}
\label{SDE}
\begin{cases}
&dX_t =  \sigma(t,X_t)\,dW_t + b(t,X_t)\,dt,\\
&X_0 = x_0
\end{cases}
\end{equation}
with $x_0 \in \R$ and $\sigma, b: [0,T] \times \R \to \R$. 
We assume that $\sigma, b \in C^{0,1}([0,T] \times \R)$ 
and for  
$f \in \lbrace \sigma,b \rbrace$ 
 there exist constants $C_T$ and $\alpha \geq \half$ such that
 
\begin{enum_i}
\item $|f(t,x)| \leq C_T(1+|x|) $,
\item $ |f(t,x)-f(t,y)| \leq C_T|x-y| $,
\item $ |f(t,x)-f(s,x)| \leq C_T(1+|x|)|t-s|^\alpha $.
\end{enum_i}
Assumptions (i)-(iii) imply the existence of a unique adapted strong solution $X$ of the SDE \eqref{SDE}, see e.g. \cite[p. 289]{KS}. For estimates concerning the Milstein scheme we assume in addition that 
\begin{enumerate}
 \item[(iv)] the state derivatives $\frac{\partial \sigma}{\partial x}$ and $\frac{\partial b}{\partial x}$ satisfy the condition (ii). 
\end{enumerate}
Moreover, we assume that
\begin{enumerate}
\item[(v)] $X_T$ has a bounded density.
\end{enumerate}

\begin{remark}
Assumption (v) is satisfied (see \cite[p. 263]{F}), if we assume that $\sigma,b \in C^\infty_b([0,T] \times \R)$ and $\sigma$ satisfies the uniform ellipticity condition, i.e. there exists a constant $\beta$ such that
$$ \sigma(t,x) \geq \beta > 0 \text{ for all } (t,x) \in [0,T] \times \R .$$
Another sufficient condition is given by Caballero et al. in \cite[Theorem 2]{CFN}. They assume that $\sigma$ and $b$ are $C^2$ in $x$, the second derivatives have polynomial growth, the functions $|\sigma(0,x)|$, $|\sigma_x(t,x)|$, $|b(0,x)|$ and $|b_x(t,x)|$ are bounded, and
$$\E \left( \left| \int_0^t \sigma(s,X_s)^2 \, ds \right|^{-p_0/2} \right) < \infty $$
for some $p_0 > 2$ and for all $t \in (0,T]$. Then there exists a continuous density $f_{X_t}$ of $X_t$ such that for all $p>1$ 
$$ f_{X_t}(x) \leq C_p \norm{ \left( \int_0^t \sigma(s,X_s)^2 \, ds \right)^{-1/2} }{p} $$
for some constant $C_p > 0$.
\end{remark}

Denote by $\pi$ a partition $0=t_0 < t_1 < \ldots < t_n = T$ of the interval $[0,T]$, and let 
$$|\pi| = \max_{0 \leq i < n} |t_{i+1}-t_i|$$ 
be the mesh size of $\pi$.
Moreover, denote an approximation of $X$ corresponding to $\pi$ by $X^{\pi}$.
Two such approximations are the well known Euler and Milstein schemes.

\begin{defn}[Euler scheme]
Let $X^E$ be the Euler scheme relative to $\pi$, i.e. 
$X^E_0=x_0$, and for $i=0, \dots , n-1$,
$$X^E_{t_{i+1}}=X^E_{t_i} + b(t_i,X^E_{t_i})(t_{i+1}-t_i)
+\sigma(t_i,X^E_{t_i})(W_{t_{i+1}}-W_{t_i}). $$
Given the values at the partition points, we also define the Euler scheme in continuous time by setting
$$X_t^E = X_{t_k}^E + \sigma(t_k,X_{t_k}^E)(W_t-W_{t_k})+b(t_k,X_{t_k}^E)(t-t_k)$$
for $t \in (t_k, t_{k+1})$.
This can be written using the integral representation
\begin{equation}
\label{Eulercontinuousscheme}
\begin{aligned}
X_t^E = x_0 &+ \int_0^t \sum_{k=0}^{n-1} \sigma(t_k,X_{t_k}^E) \chi_{[t_k,t_{k+1})}(u)\, dW_u \\
&+ \int_0^t \sum_{k=0}^{n-1} b(t_k,X_{t_k}^E) \chi_{[t_k,t_{k+1})}(u)\, du
\end{aligned}
\end{equation}
for $t \in [0,T]$ a.s., which coincides a.s. with the discrete scheme at the partition points. In this paper we need the continuous time version only for the Euler scheme.
\end{defn}

\begin{defn}[Milstein scheme]
Let $X^M$ be the Milstein scheme relative to $\pi$, i.e.
$X^M_0=x_0$, and and for $i=0, \dots , n-1$,
$$
\begin{aligned}
 X^M_{t_{i+1}}=& X^M_{t_{i}} + b(t_i,X^M_{t_{i}})(t_{i+1}-t_i) +\sigma(t_i,X^M_{t_{i}})(W_{t_{i+1}}-W_{t_i}) \\
 &+\half \sigma(t_i, X^M_{t_{i}})\sigma_x'(t_i,X^M_{t_{i}})((W_{t_{i+1}}-W_{t_i})^2-(t_{i+1}-t_i)). 
\end{aligned}
$$
\end{defn}
We remind that the Euler and Milstein approximations of $X_T$, denoted $X_T^E$ and $X_T^M$, always depend on the corresponding partition $\pi$. This is omitted from the notation for simplicity.

These assumptions hold throughout the paper.

\section{Indicator Functions}
\label{SectionErrorIndic}

\subsection{General Approximation}
\label{SubsectionGeneralApproximation}

Suppose that we have a probability space $(\Omega,\mathcal{F},\Prob)$ and two random variables $X, \hat X:\Omega \to \R$. Consider $\hat X$ to be an approximation of $X$ in the $L_p$-norm. First we find an estimate for the quantity $\E|g(X)-g(\hat X)|$, where $g$ is an indicator function.

\begin{defn}
\label{DefinitionNonIncreasingArrangement}
Recall the non-increasing rearrangement of a random variable $X$, defined by $X^*:[0,1] \to \R \cup \lbrace +\infty,-\infty \rbrace$,
$$X^*(s) := \inf\{c \in \R: \Prob(X>c)\leq s\}. $$
Here we use the convention that $\inf \emptyset = \infty$. 
\end{defn}

\begin{remark}
Definition \ref{DefinitionNonIncreasingArrangement} is slightly different from the standard non-increasing rearrangement as defined e.g. in \cite{BS}, where the absolute value of the function $X$ is taken. However, by analoguous arguments we can show the following properties:
\begin{enum_i}
 \item $X^*(1)=-\infty$, $X^*(0)=\infty$ if $X$ is not essentially bounded and $X^*(s) \in \R$ for $s \in (0,1)$,
 \item  $X^*$ is right-continuous, 
 \item $X^*$ has the same distribution as $X$ with respect to the Lebesgue measure on $[0,1]$,
\end{enum_i}
\end{remark}

\begin{defn}
\label{DefinitionDX}
 Denote the minimal slope of the function $X^*$ from the level $K$ by $d_X : \R \to [0,\infty)$,
$$d_X(K) := \inf_{\substack{s \in [0,1]\\ s \neq \alpha(K)}} \bigg\lbrace \frac{|X^*(s)-K|}{|s-\alpha(K)|} \bigg\rbrace,$$
where
$$ 
\alpha(K)= \Prob(X\geq K). 
$$
\end{defn}

\begin{thm}
\label{IndicatorError}
Suppose that $X$ is a random variable. Then the following assertions hold:
\begin{enum_i}
\item If $X$ has a bounded density $f_X$, then for all $K \in \R$, all random variables $\hat X$ and all $0 < p < \infty$ we have
$$ \E|\chi_{[K, \infty)}(X)-\chi_{[K, \infty)}(\hat X)| \leq 3D_X(K)^\frac{p}{p+1} \norm{X-\hat X}{p}^\frac{p}{p+1},$$
where
$$D_X(K) := \frac{1}{d_X(K)} \in (0, \sup f_X]. $$
Moreover, the power $\frac{p}{p+1}$ of the $L_p$-norm is optimal, i.e. if 
\begin{equation}
\label{indicthmgeneralformula}
\E|\chi_{[K, \infty)}(X)-\chi_{[K, \infty)}(\hat X)| \leq C(X,K,p) \norm{X-\hat X}{p}^\frac{p}{p+1}
\end{equation}
for all random variables $\hat X$, then the power $\frac{p}{p+1}$ can not be replaced by a power $q$ such that $\frac{p}{p+1} < q < \infty$.
\item If there exists  $p_0 > 0$ such that the formula (\ref{indicthmgeneralformula}) holds for all $p_0 \leq p < \infty$, all $K \in \R$ and all random variables $\hat X$, and there exists $B_X > 0$ such that $C(X,K,p) \leq B_X$, then $X$ has a bounded density.

\end{enum_i}
\end{thm}

\begin{proof}

Let us first show (i). Fix $K \in \R$ and $0<p<\infty$, and let $\hat X$ be a random variable such that 
$$ \E|\chi_{[K, \infty)}(X)-\chi_{[K, \infty)}(\hat X)|=\varepsilon$$
for some $\varepsilon \in (0,1]$.
 Define $\varepsilon_1 := \Prob(X \geq K,\ \hat X < K)$ and $\varepsilon_2 := \Prob(X < K,\ \hat X \geq K)$, so that $\varepsilon = \varepsilon_1 + \varepsilon_2$.  
Denote by $\alpha$ the number $\alpha(K)$ introduced in Definition \ref{DefinitionDX} 
and notice that $\alpha - \varepsilon_1 \geq 0$ and $\alpha + \varepsilon_2 \leq 1$. Now
$$ 
\begin{aligned}
\E|X-\hat X|^p &\geq \E|X-\hat X|^p\chi_{\lbrace X \geq K, \hat X< K \rbrace \cup \lbrace X < K, \hat X \geq K \rbrace } \\
&\geq \E|X-K|^p\chi_{\lbrace X \geq K, \hat X< K \rbrace \cup \lbrace X < K, \hat X \geq K \rbrace } \\
&= \E|X-K|^p\chi_{\lbrace X \geq K, \hat X< K \rbrace} 
+  \E|X-K|^p\chi_{ \lbrace X < K, \hat X \geq K \rbrace }.
\end{aligned}
$$
Since $X$ has a bounded density, we can find a number $c_0 \in [K,\infty]$ such that $\Prob(K \leq X < c_0) = \varepsilon_1$, thus also $|\lbrace K \leq X^* < c_0\rbrace |=\varepsilon_1$. Note that $c_0$ may not be unique. But $\lbrace K \leq X < c_0 \rbrace$ is a set of probability $\varepsilon_1$ where $\E|X-K|^p\chi_A$ is minimized over all $A \subset \lbrace X \geq K \rbrace$ with $\Prob(A)=\varepsilon_1$, which implies that
$$
\begin{aligned}
&\E|X-K|^p\chi_{\lbrace X \geq K, \hat X< K \rbrace} 
\geq \E|X-K|^p\chi_{[K, c_0) }(X) \\
&= \int_{[0,1]}|X^*(s)-K|^p \chi_{[K, c_0)}(X^*(s)) \, ds  
=\int_{\alpha - \varepsilon_1}^\alpha|X^*(s)-K|^p \, ds \\
&\geq \int_0^{\varepsilon_1} |d_X(K)s|^p \, ds 
= \frac{d_X(K)^p\varepsilon_1^{p+1}}{p+1}
\end{aligned}
$$
and by similar arguments
\begin{eqnarray*}
\E|X-K|^p\chi_{ \lbrace X < K, \hat X \geq K \rbrace }
\geq\int_\alpha^{\alpha + \varepsilon_2}|X^*(s)-K|^p \, ds \geq \frac{d_X(K)^p\varepsilon_2^{p+1}}{p+1}.
\end{eqnarray*}
Thus 
\begin{equation}
\label{XminusYlowerbound}
\E|X-\hat X|^p 
\geq \frac{d_X(K)^p(\varepsilon_1^{p+1}+\varepsilon_2^{p+1})}{p+1}
\geq \frac{d_X(K)^p\varepsilon^{p+1}}{2^p(p+1)}.
\end{equation}
Now the equation \eqref{XminusYlowerbound}  gives 
\begin{eqnarray*}
&&\E|\chi_{[K, \infty)}(X)-\chi_{[K, \infty)}(\hat X)| \\
&& \ \ \leq 2^\frac{p}{p+1}(p+1)^\frac{1}{p+1} \left( \frac{1}{d_X(K) }\right)^\frac{p}{p+1}\left(\E|X-\hat X|^p\right)^{\frac{1}{p+1}}.
\end{eqnarray*}
By elementary computations we can show that 
$$2^\frac{p}{p+1}(p+1)^\frac{1}{p+1} \leq 2e^\frac{1}{2e} \leq 3,$$ 
and keeping in mind the definition of $D_X$ we can write
$$ \E|\chi_{[K, \infty)}(X)-\chi_{[K, \infty)}(\hat X)| \leq 3 D_X(K)^\frac{p}{p+1}\norm{X-\hat X}{p}^\frac{p}{p+1}.$$
Using the definition of $X^*$ and the boundedness assumption for the density of $X$ we see that 
$1/d_X(K) \leq \sup f_X. $

Moreover, the power $\frac{p}{p+1}$ of $\norm{X-\hat X}{p}$ is sharp. To see this, we construct an example where the lower bound given by equation (\ref{XminusYlowerbound}) is achieved. Suppose that $\Omega = [0,1]$ is equipped with the Lebesgue measure, $K=\half$ and $\varepsilon < 1$. If we take $X(\omega)=\omega$, then $X$ has a bounded density and $d_X(\half)=1$. Now define

$$\hat X=
\begin{cases}
X, &\text{ if } \omega \in [0,\half - \frac{\varepsilon}{2}) \cup (\half + \frac{\varepsilon}{2}, 1], \\
X+\frac{\varepsilon}{2}, &\text{ if } \omega \in [\half-\frac{\varepsilon}{2},\half],\\
X-\frac{\varepsilon}{2}, &\text{ if } \omega \in (\half, \half+\frac{\varepsilon}{2}].
\end{cases}
 $$
Then
\begin{eqnarray*}
&\E|X-\hat X|^p&=\E\left|\frac{\varepsilon}{2}\right|^p\chi_{[\half - \frac{\varepsilon}{2},\half + \frac{\varepsilon}{2}]}(X) 
= \frac{\varepsilon^{p+1}}{2^p}
 = \frac{d_X(\half)^p\varepsilon^{p+1}}{2^p},
\end{eqnarray*}
which coincides with the lower bound in equation (\ref{XminusYlowerbound}) up to the constant. Hence the power $\frac{p}{p+1}$ of $\norm{X-\hat X}{p}$ can not be increased in the assertion (i).

Now we verify (ii). Let $\delta>0$ and choose $\hat X=X-\delta$. Then
$$
\begin{aligned}
\E&|\chi_{[K, \infty)}(X)-\chi_{[K, \infty)}(\hat X)| \\
&=\Prob(X \geq K,\ X-\delta < K) + \Prob(X < K,\ X-\delta \geq K) \\
&=\Prob(K \leq X < K+\delta),
\end{aligned}
$$
so that by assumption we get, for $p > p_0$, that
$$\Prob(K \leq X < K+\delta) \leq B_X \left(\E|\delta|^p\right)^{\frac{1}{p+1}} \leq B_X\delta^{\frac{p}{p+1}}. $$
We let $p$ go to infinity and conclude that
$$\Prob(K \leq X < K+\delta) \leq B_X\delta.$$
Let $N \subset \R$ be a null set with respect to the Lebesgue measure 
and let $\varepsilon > 0$. 
Since the Lebesgue outer measure of $N$ is also zero, we find a sequence 
$(I_j)$ of open intervals such that $N\subset \bigcup I_j$ and $\sum|I_j| \leq \varepsilon$. Let $\mathcal{L}_X$ be the law of $X$. Then we have
$$ \mathcal{L}_X((a,b)) \leq \mathcal{L}_X([a,b)) \leq B_X|b-a| $$
and
$$\mathcal{L}_X(N) \leq \mathcal{L}_X\left(\bigcup_j I_j\right) \leq \sum_j \mathcal{L}_X(I_j) \leq B_X\sum_j |I_j| \leq B_X\varepsilon. $$
This implies that $\mathcal{L}_X(N)=0$, so $\mathcal{L}_X$ is absolutely continuous with respect to the Lebesgue measure. By the Radon-Nikodym theorem there exists a measurable function $f:\R \to [0,\infty)$ such that
$$\mathcal{L}_X(M)=\int_M f(x)\, dx $$
for all measurable $M \subseteq \R$. Moreover, $f$ is integrable since $\mathcal{L}_X(\R)=1$. Define a function $\Phi:\R \to [0,1]$ such that
$$ \Phi(t)=\int_{-\infty}^t f(x)\, dx. $$ 
Then by  \cite[Thm. 8.17]{WR} 
we have that $\Phi'(t)=f(t)$ a.e. in $\R$.
On the other hand, we have that
$$ \Phi'(t) = \lim_{h\to 0} \frac{\Phi(t+h)-\Phi(t)}{h} \leq \lim_{h\to 0} \frac{B_Xh}{h} = B_X \text{  a.e. in } \R,$$
because $\Phi(t+h)-\Phi(t)=\mathcal{L}_X((t,t+h))$.
Therefore we conclude that $f(t) \leq B_X$ a.e. in $\R$.
\end{proof}


\begin{remark}
\label{RemarkKincludedornot}
By considering complements of the intervals in the indicator functions and the random variables $-X$ and $-\hat X$, we have corresponding results for the functions $\chi_{(K, \infty)}$, $\chi_{(-\infty, K]}$ and $\chi_{(-\infty, K)}$. 
\end{remark}

As an immediate consequence of Theorem \ref{IndicatorError}, we can derive

\begin{cor}
\label{IndicatorRateThmGeneral}
Let $X$ be the solution of the equation \eqref{SDE}, $K \in \R$ and $0 < p < \infty $. Let $X_T$ have a bounded density and suppose that $X_T^\pi$ is an approximation of $X_T$ such that
$$  \norm{X_T-X_T^\pi}{p} \leq C_p\, \mesh^\theta$$
for some $\theta > 0$ and some constant $C_p \geq 0$. Then for all $K \in \R$ we have
$$ \E|\chi_{[K, \infty)}(X_T)-\chi_{[K, \infty)}(X_T^\pi)| \leq 3D_{X_T}(K)^\frac{p}{p+1} C_p^{\frac{p}{p+1} } \,\mesh^\frac{\theta p}{p+1} .$$
\end{cor}

\subsection{Euler and Milstein Schemes}

Now we can apply the results of Section \ref{SubsectionGeneralApproximation} to the Euler and Milstein schemes:

\begin{thm}
\label{IndicatorRateThmEulerMilstein}
For any $0 < \varepsilon < 1/2$ there exists a constant $C_\varepsilon>0$ such that for all $K \in \R$ we have that 
$$\E|\chi_{[K, \infty)}(X_T)-\chi_{[K, \infty)}(X_T^E)| \leq (D_{X_T}(K) \vee \sqrt{D_{X_T}(K)}) C_\varepsilon \mesh^{\half - \varepsilon}$$ 
and for any $0 < \varepsilon < 1$ there exists a constant $C'_\varepsilon>0$ such that for all $K \in \R$ we have that
$$\E|\chi_{[K, \infty)}(X_T)-\chi_{[K, \infty)}(X_T^M)| \leq (D_{X_T}(K) \vee \sqrt{D_{X_T}(K)})C'_\varepsilon \mesh^{1 - \varepsilon}.$$
\end{thm}

\begin{proof}
Let $1 \leq p < \infty$. Then for the Euler scheme we have by Theorem \ref{Eulerthm} in the Appendix that
$$ \norm{X_T-X^E_T}{p} \leq C_p \mesh^\half,$$
i.e. the assumption of Corollary \ref{IndicatorRateThmGeneral} is satisfied with $\theta = \half$.
Thus
\begin{equation}
\label{eqindicatorEuler}
 \E|\chi_{[K, \infty)}(X_T)-\chi_{[K, \infty)}(X^E_T)| \leq 3D_{X_T}(K)^\frac{p}{p+1} C_p^\frac{p}{p+1} \mesh^\frac{p}{2(p+1)} .
\end{equation}
Similarly for the Milstein scheme we have by \cite[Proposition 1, p. 140]{M1} that
$$ \norm{X_T-X^M_T}{p} \leq C_p' \mesh ,$$
which gives the assumption of Corollary \ref{IndicatorRateThmGeneral} with $\theta = 1$, and therefore
\begin{equation}
\label{eqindicatorMilstein}
 \E|\chi_{[K, \infty)}(X_T)-\chi_{[K, \infty)}(X^M_T)| \leq 3D_{X_T}(K)^\frac{p}{p+1} (C_p')^\frac{p}{p+1} \mesh^\frac{p}{p+1} .
\end{equation}

The claim follows in both cases by choosing $p$ such that  
$p = (\theta - \varepsilon)/\varepsilon	 $,
where $0<\varepsilon < \theta$,
and noticing that for any $a>0$ we have $a^\frac{p}{p+1} \leq a \vee \sqrt{a}$. The constant $3$ and the constants coming from the approximation schemes are included in $C_\varepsilon$ or $C_\varepsilon'$, which now depend on $\varepsilon$ through the choice of $p$.
\end{proof}

Since we have information about the constant $C_p$ in Theorem \ref{Eulerthm}, i.e. $C_p=e^{Mp^2}$, we can write an extended version of Theorem \ref{IndicatorRateThmEulerMilstein} for the Euler scheme:

\begin{thm}
\label{IndicatorRateThmEulerExt}
Let $K\in \R$. Then there exists $m \in (0,1)$ such that for $\mesh  < m$ we have
$$\E|\chi_{[K, \infty)}(X_T)-\chi_{[K, \infty)}(X_T^E)| \leq (D_{X_T}(K) \vee \sqrt{D_{X_T}(K)})\mesh^{\half-\frac{2+M}{(-\log \mesh)^{1/3}}},$$
where the constant $M=M(x_0, T, C_T) \in (0,\infty)$ is taken from Theorem \ref{Eulerthm}.
\end{thm}

\begin{proof}
By Corollary \ref{IndicatorRateThmGeneral} and Theorem \ref{Eulerthm}, using
$a^\frac{p}{p+1} \leq a \vee \sqrt{a}$ for $a>0$ and $p \geq 1$, we get
\begin{equation}
\label{Egapprox1}
\begin{aligned}
 \E|\chi_{[K, \infty)}(X_T)&-\chi_{[K, \infty)}(X_T^E)| \leq 3D_{X_T}(K)^\frac{p}{p+1} C_p^{\frac{p}{p+1} } \mesh^\frac{p}{2(p+1)}   \\
 &\leq (D_{X_T}(K) \vee \sqrt{D_{X_T}(K)})3e^{Mp^2 \cdot \frac{p}{p+1}} \mesh^\frac{p}{2(p+1)}  \\
 &\leq (D_{X_T}(K) \vee \sqrt{D_{X_T}(K)})3e^{Mp^2}  \mesh^\frac{p}{2(p+1)}
\end{aligned}
\end{equation}
for all $p \geq 1$. Now choose $p$ such that
$$ 4p(p+1)^2 = -\log \mesh $$
for $\mesh \leq m$ with $m = e^{-16}$.
This gives
$ p^3 \leq -\log \mesh $ 
and $p^2 \leq  (-\log \mesh)^{2/3}$.
Thus we have

$$ e^{Mp^2} \leq e^{M(-\log \mesh)^{2/3}} = \mesh^{-M(-\log \mesh)^{-1/3}} $$
and
$$
\begin{aligned}
\frac{1}{2(p+1)} = \sqrt{\frac{p}{-\log \mesh}} \leq \sqrt{(-\log \mesh)^{1/3-1}}=(-\log \mesh)^{-1/3}.
\end{aligned}
$$
Using these we get
$$ 3e^{Mp^2}  \mesh^\frac{p}{2(p+1)}
= 3e^{Mp^2} \mesh^{\half-\frac{1}{2(p+1)}} 
\leq 3\mesh^{\half-\frac{1+M}{(-\log \mesh)^{1/3}}}
\leq \mesh^{\half-\frac{2+M}{(-\log \mesh)^{1/3}}},
$$
where in the last step we used the inequality $$ 3\mesh^\frac{1}{(-\log \mesh)^{1/3}} \leq 1 $$
for $|\pi|<m$.
\comment{
Important note:
Here must be $$ 3\mesh^\frac{1}{(-\log \mesh)^{1/3}} \leq 1, $$
which implies 
$$ \mesh^{-1} \geq e^{(\log 3)^\frac{3}{2}} , $$
so there's no need to decrease $m$, which was $e^{-16}$. 
}
Now we come back to equation \eqref{Egapprox1} and conclude that
\begin{eqnarray*}
&& \E|\chi_{[K, \infty)}(X_T)-\chi_{[K, \infty)}(X_T^E)| \\
&&\leq  (D_{X_T}(K) \vee \sqrt{D_{X_T}(K)})\mesh^{\half-\frac{2+M}{(-\log \mesh)^{1/3}}}.
\end{eqnarray*}
\end{proof}

\section{Functions of Bounded Variation}
\label{SectionErrorBV}

From Theorem \ref{IndicatorError} we deduce the same error for functions of bounded variation, up to a constant. Let us first recall the definitions of the spaces $BV$ and $NBV$. 

\begin{defn}
Let $$T_f(x) := \sup \sum_{j=1}^N |f(x_j)-f(x_{j-1})|,$$ where the supremum is taken over $N$ and all partitions $-\infty < x_0 < x_1 < \ldots < x_N = x$, be the total variation function of $f$. Then we say that $f$ is a function of bounded variation, $f \in BV$, if $$ V(f) := \lim_{x \to \infty} T_f(x) $$
is finite, and call $V(f)$ the (total) variation of $f$. 
\end{defn}

\begin{defn}
Let $NBV$ be the set of functions $f \in BV$ such that $f$ is left-continuous and $f(x) \to 0$ as $x \to -\infty$.
\end{defn}

\subsection{General Approximation}

\begin{thm}
\label{BVError}
Suppose that $X$ and $\hat X$ are random variables and $X$ has a bounded density. 
If $g \in BV$ and $1 \leq p < \infty $, then for any $1 \leq q < \infty$ we have
$$\E|g(X)-g(\hat X)|^p \leq  3^{p+1} \left(\sup f_X \right)^\frac{q}{q+1} V(g)^p \norm{X-\hat X}{q}^\frac{q}{q+1}. $$ 
\end{thm}

\begin{proof}
First we show the result for functions $g \in NBV$. By \cite[Thm. 8.14]{WR} there is a unique signed measure $\mu$ such that 
$$g(x)=\mu((-\infty,x)) \text{ and } |\mu|((-\infty,x)) = T_g(x),$$ where $|\mu|$ is the total variation measure of $\mu$. We consider the Jordan decomposition of $\mu$ , i.e. $\mu=\mu_1-\mu_2$, where $\mu_1=\half(|\mu|+\mu)$ and $\mu_2=\half(|\mu|-\mu)$ are positive measures. Then  $|\mu|=\mu_1+\mu_2$, and all three measures $|\mu|$, $\mu_1$ and $\mu_2$ are finite since $|\mu|(\R)=V(g) < \infty$.
Thus we get
$$ 
g(x)=\mu((-\infty,x))=\intR \chi_{(-\infty,x)}(z)\, d\mu(z)=\intR \chi_{(z,\infty)}(x)\, d\mu(z).
$$
Now by Theorem \ref{IndicatorError} 
and Remark \ref{RemarkKincludedornot} we get

\begin{eqnarray*}
\norm{g(X)-g(\hat X)}{p}&=&
\norm{\intR \chi_{(z,\infty)}(X)\, d\mu(z)-\intR \chi_{(z,\infty)}(\hat X)\, d\mu(z)}{p} \\
&=&\norm{\intR \left[ \chi_{(z,\infty)}(X) - \chi_{(z,\infty)}(\hat X) \right] \, d\mu(z)}{p} \\
&\leq& \norm{\intR \left|\chi_{(z,\infty)}(X) - \chi_{(z,\infty)}(\hat X)\right|\, d|\mu|(z)}{p} \\
&\leq& \intR \norm{\chi_{(z,\infty)}(X) - \chi_{(z,\infty)}(\hat X)}{p}\, d|\mu|(z) \\
&\leq& 3^\frac{1}{p}  (\sup f_X)^\frac{q}{p(q+1)} V(g) \norm{X-\hat X}{q}^\frac{q}{p(q+1)},
\end{eqnarray*}
which completes the proof for functions in $NBV$.

Next, let $g$ be an arbitrary function in $BV$. By \cite[Thm. 8.13]{WR}, there exists a unique function $\tilde g \in NBV$ and a unique constant $c \in \R$ such that $g(x) = \tilde g(x) + c$ at all points of continuity of $g$, with $V(\tilde g) \leq V(g)$. Also by \cite{WR} we know that 
$g$ can have only countably many points of discontinuity, so define $\cup_{j=1}^\infty \lbrace a_j \rbrace$ to be the set of these points and let $\lambda_j:=g(a_j)-\tilde g(a_j) - c$. 
Then we can write
$$ g(x)= \tilde g(x) + c + \Delta(x), $$
where  
$$ \Delta(x):= \sum_{j=1}^\infty \lambda_j \chi_{\lbrace a_j \rbrace}(x)  = \sum_{j=1}^\infty \lambda_j \left(\chi_{(-\infty,a_j]}(x)- \chi_{(-\infty,a_j)}(x) \right).$$
We define a measure 
$$\nu=\sum_{j=1}^\infty \lambda_j \delta_{a_j},$$
where $\delta_a$ is the Dirac measure in $a$. 
Again by \cite{WR} we know that $g(a_j-)$ exists, so we have $\tilde g(a_j)+c=g(a_j-)$ and
$$|\nu|(\R)
=\sum_{j=1}^\infty |\lambda_j| = \sum_{j=1}^\infty |g(a_j)-g(a_j-)| \leq V(g). $$
Now we can write
$$ \Delta(x)=\intR  \chi_{(-\infty,z]}(x)- \chi_{(-\infty,z)}(x) \, d\nu(z).$$
and compute, similarly as in the NBV case, that
\begin{eqnarray*}
&&\norm{\Delta(X)-\Delta(\hat X)}{p} 
\leq  \norm{\intR  |\chi_{(-\infty,z]}(X)- \chi_{(-\infty,z]}(\hat X)| \, d|\nu|(z)}{p} \\ 
&&\ \ +\norm{\intR  |\chi_{(-\infty,z)}(X)- \chi_{(-\infty,z)}(\hat X)| \, d|\nu|(z)}{p} \\
&&  \leq \intR  \norm{\chi_{(-\infty,z]}(X)- \chi_{(-\infty,z]}(\hat X)}{p} \, d|\nu|(z) \\
&&\ \   +\intR  \norm{\chi_{(-\infty,z)}(X)- \chi_{(-\infty,z)}(\hat X)}{p}  \, d|\nu|(z) \\
&&\leq  2 \cdot 3^\frac{1}{p} (\sup f_X)^\frac{q}{p(q+1)} V(g) \norm{X-\hat X}{q}^\frac{q}{p(q+1)}.
\end{eqnarray*}
This, combined with the NBV result, implies that
$$ 
\begin{aligned}
&\norm{g(X)-g(\hat X)}{p} 
= \norm{\tilde g(X)- \tilde g(\hat X) + \Delta(X) - \Delta(\hat X)}{p}  \\
&\leq \norm{\tilde g(X)- \tilde g(\hat X)}{p} + \norm{\Delta(X) - \Delta(\hat X)}{p} \\
&\leq 3 \cdot 3^\frac{1}{p}  (\sup f_X)^\frac{q}{p(q+1)} V(g) \norm{X-\hat X}{q}^\frac{q}{p(q+1)}, 
\end{aligned}
$$
which gives the statement.
\end{proof}

As in Corollary \ref{IndicatorRateThmGeneral} for indicator functions,
we can now write an analoguous statement for functions of bounded variation:

\begin{cor}
\label{BVRateThmGeneral}
Let $X$ be the solution of the equation \eqref{SDE}, $1 \leq p < \infty $ and $g \in BV$. Suppose that $X_T$ has a bounded density, $1\leq q < \infty$ and $X_T^\pi$ is an approximation of $X_T$ such that 
$$  \norm{X_T-X_T^\pi}{q} \leq C_q \mesh^\theta$$
for some $\theta > 0$ and some constant $C_q \geq 0$. Then
$$ \E|g(X_T)-g(X_T^\pi)|^p \leq 3^{p+1} \left(\sup f_{X_T} \right)^\frac{q}{q+1} V(g)^p C_q^\frac{q}{q+1} \mesh^\frac{\theta q}{q+1} .$$
\end{cor}

\subsection{Euler and Milstein Schemes}


\begin{thm}
\label{BVRateThmEulerMilstein}
Let $g \in BV$ and $1 \leq p < \infty$. Then we have for $0 < \varepsilon < 1/2$ that
$$\E|g(X_T)-g(X_T^E)|^p \leq 3^p(\sup f_{X_T} \vee \sqrt{\sup f_{X_T}})V(g)^p C_\varepsilon \mesh^{\half - \varepsilon}$$
and for $0 < \varepsilon < 1$ that
$$\E|g(X_T)-g(X_T^M)|^p \leq 3^p (\sup f_{X_T} \vee \sqrt{\sup f_{X_T}})V(g)^p C_\varepsilon' \mesh^{1 - \varepsilon},$$
where $C_\varepsilon$ and $C_\varepsilon'$ depend on $\varepsilon$ and the constants of the corresponding schemes.
\end{thm}

\begin{proof}
The statement follows from Corollary \ref{BVRateThmGeneral} with arguments similar to the proof of Theorem \ref{IndicatorRateThmEulerMilstein}.
\end{proof}

For the Euler scheme we can write an extended version corresponding to Theorem \ref{IndicatorRateThmEulerExt}:

\begin{thm}
\label{BVRateThmEulerExt}
Let $1 \leq p < \infty $ and $g \in BV$. Then there exists $m > 0$ such that for $\mesh < m$ we have
$$\E|g(X_T)-g(X_T^E)|^p \leq 3^p (\sup f_{X_T} \vee \sqrt{\sup f_{X_T}})V(g)^p \mesh^{\half-\frac{2+ M}{(-\log \mesh)^{1/3}}},$$ 
where $M$ is the constant in Theorem \ref{Eulerthm}.
\end{thm}
\begin{proof}
By Theorem \ref{Eulerthm} and Corollary \ref{BVRateThmGeneral} we get for $1 \leq q < \infty$ that 
$$ 
\begin{aligned}
\E|g(X_T)-g(X^E_T)|^p 
&\leq 3^{p+1} (\sup f_{X_T}\!\! \vee\! \sqrt{\sup f_{X_T}})V(g)^p e^{M q^2 \cdot \frac{q}{q+1}} \mesh^\frac{q}{2(q+1)}  \\
&\leq   3^{p+1} (\sup f_{X_T} \vee \sqrt{\sup f_{X_T}})V(g)^p e^{M q^2} \mesh^\frac{q}{2(q+1)}, 
\end{aligned}
$$
which by the arguments in Theorem \ref{IndicatorRateThmEulerExt} implies that
$$
\E|g(X_T)-g(X^E_T)|^p \leq 3^p (\sup f_{X_T} \vee \sqrt{\sup f_{X_T}})V(g)^p \mesh^{\half-\frac{2+ M}{(-\log \mesh)^{1/3}}}.
$$
\end{proof}

\section{Extension}
\label{SectionExtension}

Next we extend the result into a function class, to be called $\G_{p,\varphi}$, that contains e.g. all polynomials. The main result is given in Theorem \ref{Gpthm} and the class $\G_{p,\varphi}$ is analyzed in Section \ref{SectionGpAnalysis}.

\begin{defn}[Bump function]
Let $\varphi: \R \to \R$ be a function such that 
$ 0 < \varphi(z) \leq 1$ for all $z \in \R$, $\varphi$ is increasing in $(-\infty,0]$ and decreasing in $(0,\infty)$, and $$\lim_{|z|\to \infty}\varphi(z)=0.$$ Then $\varphi$ is called a bump function.
\end{defn}

\begin{defn}[Class $\G_{p,\varphi}$]
\label{Gpdef}
Fix $p \in [1,\infty)$ and let $\varphi: \R \to \R$ be a given bump
function. 
Let $\mathcal{M}$ be the set of all signed measures $\mu$ on $(\R,\mathcal{B}(\R))$ such that $|\mu|$ is $\sigma$-finite. Define 
$$\mathcal{M}_{p,\varphi} = \lbrace \mu \in \! \mathcal{M} :  \varphi \in 
L_1^{loc}(\R,|\mu|) \cap L_{1+\frac{1}{p}}(\R,|\mu|) \rbrace. $$
Then for any $\mu \in \mathcal{M}_{p,\varphi}$ define a function related to $\mu$ by
$$ g^\mu(x)= 
\begin{cases}
\int_{(0,x]}\varphi \, d\mu ,\ \ \text{ for } x \geq 0,\\
\int_{(x,0]}\varphi \, d\mu ,\ \  \text{ for } x < 0,
\end{cases}
 $$
where $(0,0]=\emptyset$. 
Also define a set of jump functions
\begin{eqnarray*}
&\Delta_{p,\varphi} = \lbrace& \!\!\Delta_A(x)=\sum_{i=1}^\infty \lambda_i \varphi(a_i) \chi_{\lbrace a_i \rbrace }(x) : A=(a_i)_{i=1}^\infty \subset \R \text{ countable }, \\
&& a_i \neq a_j \text{ if } i \neq j,\ (\lambda_i)_{i=1}^\infty \subset \R \text{ and } \sum_{i=1}^\infty |\lambda_i| \varphi(a_i)^{1+\frac{1}{p}} < \infty \rbrace.
\end{eqnarray*}
\comment{
For a countable set $A=(a_i)_{i=1}^\infty \subset \R$ define a jump function 
$$\Delta_A(x)=\sum_{i=1}^\infty \lambda_i \chi_{\lbrace a_i \rbrace }(x),$$
where the coefficients $\lambda_i \in \R$ satisfy the condition
$$ \sum_{i=1}^\infty |\lambda_i| \varphi(a_i)^{1+\frac{1}{p}} < \infty. $$ }
Then we set
\begin{eqnarray*}
& \G_{p,\varphi} = \lbrace g:& g=c+g^\mu + \Delta_A \text{ for some } c \in \R, \\
&& \mu \in \mathcal{M}_{p,\varphi} \text{ and } \Delta_A \in \Delta_{p,\varphi} \rbrace, 
\end{eqnarray*}
where the decomposition of $g$ is unique, as we will see in Theorem \ref{decompositionunique}.
Moreover, denote the $(p,\varphi)$-variation of $g \in \G_{p,\varphi}$ by
$$V_{p,\varphi}(g) = \int_{\R} \varphi^{1+\frac{1}{p}} \, d|\mu| + \sum_{i=1}^\infty |\lambda_i| \varphi(a_i)^{1+\frac{1}{p}}. $$
\end{defn}

\begin{remark}
 The definition implies that any function $g^\mu \in \G_{p,\varphi}$ is right-continuous and $g^\mu(0)=0$. To relax these restrictions, we add to the function $g^\mu$ a constant $c$ and a function $\Delta_A$, which can be used to alter the left- or right-continuity of $g^\mu$ at the points of discontinuity or to add point discontinuities anywhere. For example, we can make $g^\mu$ left-continuous by choosing
$$\Delta_A(x)=
\begin{cases}
 g^\mu(x-)-g^\mu(x) \text{ for } x \in A, \\
0 \text{ elsewhere, }
\end{cases}
$$
where $A$ is the set of points of discontinuity of $g^\mu$. We see that $g$ can have only a countable number of jumps, because otherwise $\mu(\lbrace x \rbrace) \neq 0$ for uncountably many $x \in \R$, which is a contradiction to the $\sigma$-finiteness of the measure $|\mu|$. Moreover, both $g$ and $g^\mu$ may jump at zero; the jump height of $g^\mu$ is then given by $\varphi(0)\mu(\lbrace 0 \rbrace)$, and the jump of $g$ depends on $\varphi(0)\mu(\lbrace 0 \rbrace)$ and $\Delta_A(0)$.
\end{remark}

\begin{thm}
\label{decompositionunique}
The decomposition $g=c+g^\mu + \Delta_A$ for functions $g \in \G_{p,\varphi}$  is unique.
\end{thm}
\begin{proof}
Take $g_1,g_2 \in \G_{p,\varphi}$ such that $g_i=c_i+g^{\mu_i} + \Delta_{A_i}$, $i \in \lbrace 1,2 \rbrace$, and suppose that $g_1=g_2$. Now $A_1 \cup A_2$ is countable and $\Delta_{A_i}=0$ in $(A_1 \cup A_2)^c$. Let us take a sequence $(x_j) \subset (A_1 \cup A_2)^c$ such that $x_j \searrow 0$ as $j \to \infty$. Since $g^{\mu_i}$ is right-continuous and $g^{\mu_i}(0)=0$, we get that $g_i(x_j)=c_i + g^{\mu_i}(x_j) \to c_i$, and thus $c_1=c_2$.
This implies that for $x_0 \in (A_1 \cup A_2)^c$ we have $g^{\mu_1}(x_0)=g^{\mu_2}(x_0)$. Now let $x_0 \in A_1 \cup A_2$. Again we choose a sequence $(x_j) \subset (A_1 \cup A_2)^c$ such that $x_j \searrow x_0$ as $j \to \infty$, and by right-continuity of $g^{\mu_i}$ we get that $g^{\mu_1}(x_0)=g^{\mu_2}(x_0)$. Thus $g^{\mu_1}=g^{\mu_2}$ everywhere, and also $\Delta_{A_1}=\Delta_{A_2}$.
\end{proof}

\begin{thm}
Functions of bounded variation are a special case of functions in $\G_{p,\varphi}$, i.e. $BV \subset \G_{p,\varphi}$.
\end{thm}
\begin{proof}
Let $g \in NBV$ and let $\mu^{BV}$ be the signed measure related to $g$. At the points of continuity of $g$ we have $g=g(0+)+g^\mu$, where the measure $\mu$ is chosen such that $d\mu=d\mu^{BV}/\varphi$ on $(0,\infty)$ and $d\mu=-d\mu^{BV}/\varphi$ on $(-\infty,0]$. Here $|\mu|$ is $\sigma$-finite by the finiteness of $\mu^{BV}$ and the properties of $\varphi$. It also holds that $$\norm{\varphi}{L_{1+\frac{1}{p}}(\R,|\mu|)} \leq V(g)^\frac{p}{p+1} < \infty.$$ 
Now let $g \in BV$. Then $g=\tilde g + c$ for some $\tilde g \in NBV$ and $c \in \R$ at the points of continuity of $g$, thus satisfying $g=\tilde g(0+)+c+\tilde g^\mu$. At the points of discontinuity we correct this by choosing $\Delta_A$ such that $A$ is the set of the points where $g$ is not right-continuous, and the values $\lambda$ correspond to the jump heights of $g$. Then
$$\sum_{i=1}^\infty |\lambda_i| \varphi(a_i)^{1+\frac{1}{p}} \leq \sum_{i=1}^\infty |\lambda_i| \leq V(g).$$
\end{proof}

\subsection{General Approximation}

As before, let $X$ and $\hat X$ be random variables. We define a function $\varphi$ that connects the random variables with their tail behavior.

\begin{defn}
\label{DEFphitailestimate}
Take two strictly positive monotone functions $$\varphi_\theta^+:(0,\infty) \to (0,1]$$ and $$\varphi_\theta^-:(-\infty,0] \to (0,1]$$ for $0 < \theta < 1$ with properties
$$\varphi_\theta^+(K) \to 0 \text{ as } K \to \infty,$$
$$\varphi_\theta^-(K) \to 0 \text{ as } K \to -\infty,$$ 
$$\left[\Prob(X \geq K) \vee \Prob(\hat X \geq K) \right]^\theta \leq \varphi_\theta^+(K) \text{ for } K > 0$$
and
$$\left[\Prob(X \leq K) \vee \Prob(\hat X \leq K) \right]^\theta \leq \varphi_\theta^-(K) \text{ for } K \leq 0.$$
Then we define a bump function $\varphi_\theta^{X,\hat X}:\R \to (0,1]$ by
$$\varphi_\theta^{X,\hat X}(K) :=
\begin{cases}
 \varphi_\theta^+(K) &\text{ if } K > 0, \\
 \varphi_\theta^-(K) &\text{ if } K \leq 0. 
\end{cases}
$$
\end{defn}

Now the main result is the following convergence theorem for functions in the class $\G_{p,\varphi}$ associated with the function $\varphi_\theta^{X,\hat X}$.

\begin{thm}
\label{Gpthm}
Let $X$ and $\hat X$ be random variables such that $X$ has a bounded density. Suppose that $0<\theta<1$ and let $\varphi_\theta^{X,\hat X}$ be a function as in Definition \ref{DEFphitailestimate}. If $1\leq p < \infty$ and $g \in \G_{p,\varphi_\theta^{X,\hat X}}$, then for all $1 \leq q < \infty$ we have
$$\E|g(X)-g(\hat X)|^p \leq 3 \cdot 2^p \left(\sup f_X \right)^\frac{q(1-\theta)}{q+1} \left(V_{p,\varphi_\theta^{X,\hat X}}(g)\right)^p \norm{X-\hat X}{q}^\frac{q(1-\theta)}{q+1}. $$
\end{thm}

\begin{proof}
Let $g \in \G_{p,\varphi_\theta^{X,\hat X}}$. Then by definition $g=c+ g^\mu + \Delta_A$ and 
$$\norm{g(X)-g(\hat X)}{p} \leq \norm{g^\mu(X)-g^\mu(\hat X)}{p} + \norm{\Delta_A(X)-\Delta_A(\hat X)}{p}. $$
Now we 
can compute

\begin{eqnarray*}
g^\mu(x)\chi_{[0,\infty)}(x)
&=&\int_{(0,x]}\varphi_\theta^{X,\hat X}(z) \, d\mu(z) \ \chi_{[0,\infty)}(x) \\
&=&\int_{(0,\infty)}\chi_{[z,\infty)}(x)\varphi_\theta^{X,\hat X}(z) \, d\mu(z) 
\end{eqnarray*}
and similarly
\begin{eqnarray*}
g^\mu(x)\chi_{(-\infty,0)}(x)
&=&\int_{(x,0]}\varphi_\theta^{X,\hat X}(z) \, d\mu(z) \ \chi_{(-\infty,0)}(x) \\
&=&\int_{(-\infty,0]}\chi_{(-\infty,z)}(x)\varphi_\theta^{X,\hat X}(z) \, d\mu(z).
\end{eqnarray*}
Thus
\begin{eqnarray*}
&&\norm{g^\mu(X)-g^\mu(\hat X)}{p}  
\leq \norm{g^\mu(X)\chi_{[0,\infty)}(X)-g^\mu(\hat X)\chi_{[0,\infty)}(\hat X)}{p} \\
 &&\ + \norm{g^\mu(X)\chi_{(-\infty,0)}(X)-g^\mu(\hat X)\chi_{(-\infty,0)}(\hat X)}{p} \\
&&\leq \norm{\int_{(0,\infty)}  |\chi_{[z,\infty)}(X)-\chi_{[z,\infty)}(\hat X)|\varphi_\theta^{X,\hat X}(z) \, d|\mu|(z)}{p} \\
&&\ + \norm{\int_{(-\infty,0]}  |\chi_{(-\infty,z)}(X)-\chi_{(-\infty,z)}(\hat X)|\varphi_\theta^{X,\hat X}(z) \, d|\mu|(z)}{p} \\
&&\leq \int_{(0,\infty)}  \norm{\chi_{[z,\infty)}(X)-\chi_{[z,\infty)}(\hat X)}{p}\varphi_\theta^{X,\hat X}(z) \, d|\mu|(z) \\
&&\ + \int_{(-\infty,0]}  \norm{\chi_{(-\infty,z)}(X)-\chi_{(-\infty,z)}(\hat X)}{p}\varphi_\theta^{X,\hat X}(z) \, d|\mu|(z). 
\end{eqnarray*}
Denote by $\psi(X,\hat X)$ the error function from Theorem \ref{IndicatorError}, i.e.
$$ \psi(X,\hat X) :=  3(\sup f_X)^\frac{q}{q+1}\norm{X-\hat X}{q}^\frac{q}{q+1}$$
and notice that $a \wedge b \leq a^{1-\theta}b^\theta$ for any $a,b \geq 0$ and $0 < \theta < 1$. Since
$$
\begin{aligned}
\E|\chi_{[K, \infty)}(X)-\chi_{[K, \infty)}(\hat X)| 
&=\Prob(X \geq K, \hat X \! < K) + \Prob(X < K, \hat X \geq K)\\
&\leq 2(\Prob(X \geq K) \vee \Prob(\hat X \geq K)),
\end{aligned}
$$
it follows from Theorem \ref{IndicatorError} that, for $K > 0$,
$$
\begin{aligned}
\E|\chi_{[K,\infty)}(X)-\chi_{[K,\infty)}(\hat X)| &\leq \psi(X,\hat X) \wedge 2 \left[ \Prob(X \geq K) \vee \Prob(\hat X \geq K) \right]\\
&\leq \psi(X,\hat X)^{1-\theta}2^\theta \! \left[\Prob(X \geq K) \! \vee \Prob(\hat X  \geq K) \right]^\theta \\
&\leq 2^\theta \psi(X,\hat X)^{1-\theta}\varphi_\theta^+(K),
\end{aligned}
$$
where $0 < \theta < 1$. In a similar way we get for $K \leq 0$ that
$$
\E|\chi_{[K,\infty)}(X)-\chi_{[K,\infty)}(\hat X)| \leq 2^\theta \psi(X,\hat X)^{1-\theta}\varphi_\theta^-(K),
$$
so we can write for $K \in \R$ that
$$
\E|\chi_{[K,\infty)}(X)-\chi_{[K,\infty)}(\hat X)| \leq 2^\theta \psi(X,\hat X)^{1-\theta}\varphi_\theta^{X,\hat X}(K).
$$
This gives an estimate for $\norm{\chi_{[z,\infty)}(X)-\chi_{[z,\infty)}(\hat X)}{p}$, and the same estimate holds for $\norm{\chi_{(-\infty,z)}(X)-\chi_{(-\infty,z)}(\hat X)}{p}$ by the observation in
Remark \ref{RemarkKincludedornot}. 
Therefore
\begin{eqnarray*}
&&\norm{g^\mu(X)-g^\mu(\hat X)}{p}  
\leq \intR  2^\frac{\theta}{p} \psi(X,\hat X)^\frac{1-\theta}{p}\varphi_\theta^{X,\hat X}(z)^\frac{1}{p} \varphi_\theta^{X,\hat X}(z)\, d|\mu|(z) \\
&&\leq 2^\frac{\theta}{p} \int_\R 3^\frac{1-\theta}{p} (\sup f_X)^\frac{q(1-\theta)}{(q+1)p} \norm{X-\hat X}{q}^\frac{q(1-\theta)}{(q+1)p}\varphi_\theta^{X,\hat X}(z)^{1+\frac{1}{p}} \, d|\mu|(z) \\
&&\leq 3^\frac{1}{p} (\sup f_X)^\frac{q(1-\theta)}{(q+1)p} \int_\R \varphi_\theta^{X,\hat X}(z)^{1+\frac{1}{p}} \, d|\mu|(z) \norm{X-\hat X}{q}^\frac{q(1-\theta)}{(q+1)p}.
\end{eqnarray*}

It remains to show a similar estimate for the jump function $\Delta_A$. This can be done by the same argument as in the case of bounded variation, namely by writing
$$ \Delta_A(x)=\intR  \chi_{(-\infty,z]}(x)- \chi_{(-\infty,z)}(x) \, d\nu(z),$$
where 
$$\nu=\sum_{i=1}^\infty \lambda_i \varphi_\theta^{X,\hat X}(a_i) \delta_{a_i}$$
and $\delta_a$ is the Dirac measure in $a$. 
Then by arguments similar to the first part of the proof and Remark \ref{RemarkKincludedornot} we get
\begin{eqnarray*}
&&\norm{\Delta_A(X)-\Delta_A(\hat X)}{p} 
 \leq  \norm{\intR  |\chi_{(-\infty,z]}(X)-\chi_{(-\infty,z]}(\hat X)|\, d|\nu|(z)}{p} \\
 & &\ \  + \norm{\intR  |\chi_{(-\infty,z)}(X)-\chi_{(-\infty,z)}(\hat X)|\, d|\nu|(z)}{p} \\
&&\leq \intR  \norm{\chi_{(-\infty,z]}(X)-\chi_{(-\infty,z]}(\hat X)}{p}\, d|\nu|(z) \\
&&+ \intR  \norm{\chi_{(-\infty,z)}(X)-\chi_{(-\infty,z)}(\hat X)}{p}\, d|\nu|(z) \\
&& \leq 2 \cdot  3^\frac{1}{p} (\sup f_X)^\frac{q(1-\theta)}{(q+1)p} \int_\R \varphi_\theta^{X,\hat X}(z)^\frac{1}{p} \, d|\nu|(z) \norm{X-\hat X}{q}^\frac{q(1-\theta)}{(q+1)p} \\
&& \leq 2 \cdot  3^\frac{1}{p} (\sup f_X)^\frac{q(1-\theta)}{(q+1)p} \left(\sum_{i=1}^\infty |\lambda_i| \varphi_\theta^{X,\hat X}(a_i)^{1+\frac{1}{p}}\right) \norm{X-\hat X}{q}^\frac{q(1-\theta)}{(q+1)p},
\end{eqnarray*}
so finally we get
$$
\begin{aligned}
\norm{g(X)-g(\hat X)}{p} \leq 
  2 \cdot  3^\frac{1}{p} (\sup f_X)^\frac{q(1-\theta)}{(q+1)p} V_{p,\varphi_\theta^{X,\hat X}}(g) \norm{X-\hat X}{q}^\frac{q(1-\theta)}{(q+1)p}.
\end{aligned}
$$
\end{proof}

\section{Analysis of the Class $\G_{p, \varphi}$}
\label{SectionGpAnalysis}
We study the class $\G_{p, \varphi}$ with the underlying function $\varphi=\varphi_\theta^{X,\hat X}$. This function depends on the approximation $\hat X$, and our first task is to handle this dependence. We show in Lemma \ref{Guniformboundforphi} that we can choose the function  $\varphi_\theta^{X,\hat X}$ such that it decays faster than any polynomial, and then we prove in Theorem \ref{polynomialsinG} that with this choice, the class $\G_{p, \varphi_\theta^{X,\hat X}}$ contains all polynomials. Then we apply the results to solutions of SDEs, and collect our knowledge in the main result, Corollary \ref{ExtensionRateCorGeneral}.

\begin{lemma} Suppose that $\varphi$ and $\psi$ are bump functions.
\label{Gpordered}
\hfill
\begin{enum_i}
\item If $\varphi \leq \psi$, then $\G_{p, \psi} \subset \G_{p, \varphi}$.
\item If $g \in \G_{p, \psi}$ and $\varphi \leq \psi$, then $V_{p,\varphi}(g) \leq V_{p,\psi}(g)$.
\end{enum_i}
\end{lemma}
\begin{proof}
First we show (i). Let $g \in \G_{p, \psi}$ and $\mu_\psi$ be related to $g$, i.e. $g=c + g^{\mu_\psi} + \Delta_A$. We choose a measure $\mu_\varphi$ such that $d\mu_\varphi = (\psi/\varphi) d\mu_\psi$, which implies that $d|\mu_\varphi| = (\psi/\varphi) d|\mu_\psi|$ and $|\mu_\varphi|$ is $\sigma$-finite. 
Now we get for $x\geq 0$ that
$$g^{\mu_\psi}(x)=\int_{(0,x]} \psi(z)\, d\mu_\psi(z) = \int_{(0,x]} \varphi(z) \, d\mu_\varphi(z) $$
and similarly for $x<0$. The integrability conditions are satisfied, since
$$\int_A \varphi(z) \, d|\mu_\varphi|(z)=\int_A \psi(z) \, d|\mu_\psi|(z) < \infty $$
for all $A \subset \subset \R$ and
$$
\begin{aligned}
\intR \varphi(z)^{1+\frac{1}{p}} \, d|\mu_\varphi|(z)
&=\intR \psi(z)\varphi(z)^\frac{1}{p} \, d|\mu_\psi|(z) \\
&\leq \intR \psi(z)^{1+\frac{1}{p}} \, d|\mu_\psi|(z) < \infty. 
\end{aligned}
$$
The representation of the jump part $\Delta_A$ changes correspondingly in the change of measure, i.e. we set $$\lambda_i^\varphi=\lambda_i^\psi \cdot \frac{\psi(a_i)}{\varphi(a_i)}$$ and see that
$$\sum_{i=1}^\infty |\lambda_i^\varphi| \varphi(a_i)^{1+\frac{1}{p}} 
\leq \sum_{i=1}^\infty |\lambda_i^\psi| \psi(a_i)^{1+\frac{1}{p}} < \infty.$$
This proves the assertion (i), and (ii) follows by a similar argument.
\end{proof}

\begin{lemma}
\label{Guniformboundforphi}
Suppose that $X$ and $\hat X$ are random variables such that $X \in \bigcap_{p \in [1, \infty)} L_p$, and suppose there exists $C=(C_p)_{p \in [1,\infty)} \subset (0,\infty)$ such that $\norm{X-\hat X}{p} \leq C_p$ for all $p \in [1,\infty)$. Let $\theta \in (0,1)$. Then we can choose the function $\varphi_\theta^{X,\hat X}$ such that
$\varphi_\theta^{X,\hat X} = \varphi^X_{C,\theta}$, where the function $\varphi^X_{C,\theta}$ is a bump function that decays faster than any polynomial. 
\end{lemma}
\begin{proof}
The triangle inequality gives that $\hat X \in L_p$ and 
$$\norm{\hat X}{p} \leq \norm{X-\hat X}{p} + \norm{X}{p} \leq C_p + \norm{X}{p}. $$
Thus by Chebychev's inequality we have  for all $\lambda > 0$ that
$$\Prob(|X| \geq \lambda) \leq \frac{\E|X|^p}{\lambda^p}$$
and
$$\Prob(|\hat X| \geq \lambda) \leq  \frac{\E|\hat X|^p}{\lambda^p} \leq \frac{(C_p + \norm{X}{p})^p}{\lambda^p}. $$
So we have a polynomial tail estimate for $X$ and $\hat X$ that depends only on the constants $C_p$ of the $L_p$-estimates, not directly on $\hat X$. This implies that 
$$\left[\Prob(|X| \geq \lambda) \vee \Prob(|\hat X| \geq \lambda) \right]^\theta \leq \inf_{p \in [1, \infty)} \frac{(C_p + \norm{X}{p})^{\theta p}}{\lambda^{\theta p}} \wedge 1 =: \varphi_{C,\theta}^{X,0}(\lambda) $$ 
for $\lambda > 0$. For $\lambda < 0$ we define $\varphi_{C,\theta}^{X,0}(\lambda):=\varphi_{C,\theta}^{X,0}(|\lambda|)$, and $\varphi_{C,\theta}^{X,0}(0) := 1$.  The function $\varphi_{C,\theta}^{X,0}$ satisfies the monotonicity properties of a bump function, but is not necessarily strictly positive.
However, if we take a bump function $\psi$ and define $$\varphi^X_{C,\theta} := \varphi_{C,\theta}^{X,0} \vee \psi, $$
then $\varphi^X_{C,\theta}$ is a bump function suitable for the choice of $\varphi_\theta^{X,\hat X}$. 
Since $\varphi_{C,\theta}^{X,0}$ clearly decays faster than any polynomial and we can choose $\psi(\lambda)=e^{-|\lambda|}$, we see that $\varphi^X_{C,\theta}$ also decays faster than any polynomial.
\end{proof}

Let $\Pcal$ be the set of all polynomials from $\R$ to $\R$. Then we have the following:

\begin{thm}
\label{polynomialsinG}
Suppose that $\varphi$ is a bump function that decays faster than any polynomial. Then $\Pcal \subset \G_{p,\varphi}$ for all $p \in [1,\infty)$.
\comment{
Version with all possible assumptions:
Let $X$ be a random variable with a bounded density and assume that $X \in L^p$. Let $Y^n$ to be a sequence of random variables such that 
$$\norm{X-Y^n}{p} \leq C_p\, n^{-\gamma}$$ for some $\gamma > 0$ and $C_p > 0$, and for all $1 \leq p < \infty$.
Then $\Pcal \subset \G_{q,\varphi_\theta^{X,Y^n}}$ for all $1 \leq q< \infty$.
}
\end{thm}
\begin{proof}
Let $p \in [1,\infty)$ and suppose that $g \in \Pcal$. Then by the fundamental theorem of calculus we have for $x > 0$ that
$$g(x)=g(0) + \int_0^x g'(z) \, dz $$
and for $x \leq 0$ that
$$g(x)=g(0) - \int_x^0 g'(z) \, dz. $$
Thus by defining $c=g(0)$ and a signed measure $\mu$ such that
$$d\mu(z)=\sgn(z) \frac{g'(z)}{\varphi(z)} \,dz,$$
we have that $|\mu|$ is $\sigma$-finite and the representation $g=c+g^\mu$ holds.
\comment{
$\mu=\mu_1 + \mu_2$, where
$$d\mu_1=\frac{g'\chi_{(0,\infty)}}{\varphi} \,dz$$
and
$$d\mu_2=-\frac{g'\chi_{(-\infty,0]}}{\varphi}\,dz$$
are $\sigma$-finite signed measures,
}
Now $g'$ also has only polynomial growth, say $|g'(x)| \leq C(1+|x|^s)$ for $s \geq 1$. But $\varphi$ decays faster than any polynomial, so we have $\varphi(x) \leq \tilde C|x|^{-p(s+2)} \wedge 1$ and
\begin{eqnarray*}
&V_{p,\varphi}(g)=\int_{\R} \varphi^{1+\frac{1}{p}} \, d|\mu|=\int_{\R} \varphi^{\frac{1}{p}}(z)|g'(z)| \, dz \\ 
&\hspace{45pt} \leq C \tilde C^\frac{1}{p} \int_{\R} (|z|^{-(s+2)} \wedge 1)(1+|z|^s) \, dz
 < \infty, 
\end{eqnarray*}
which implies that $g \in \G_{p,\varphi}$.
\end{proof}

Let us now come back to the SDE \eqref{SDE} and summarize our knowledge:

\begin{cor}
\label{ExtensionRateCorGeneral}
Let $p \in [1, \infty)$. Suppose that $X$ is the solution of the equation \eqref{SDE}, $X_T$ has a bounded density, and $X_T^\pi$ is an approximation of $X_T$ such that 
$$  \norm{X_T-X_T^\pi}{p} \leq C_p \mesh^\gamma$$
for some constants $\gamma > 0$ and $C_p \geq 0$. Then for any $0 < \varepsilon < \gamma$ we have for $\theta = \frac{\varepsilon}{2\gamma - \varepsilon}$, $\varphi_{C,\theta}^{X_T}$ according to Lemma \ref{Guniformboundforphi} and 
$g \in \G_{p,\varphi_{C,\theta}^{X_T}}$ that
$$ \E|g(X_T)-g(X_T^\pi)|^p \leq 3 \cdot 2^p (\sup f_{X_T})^{1 - \frac{\varepsilon}{\gamma}} \left(V_{p,\varphi_{C,\theta}^{X_T}}(g)\right)^p C_{1/\theta}^{1 - \frac{\varepsilon}{\gamma}} \mesh^{\gamma - \varepsilon}.$$
Especially, $\Pcal \subset \G_{p,\varphi_{C,\theta}^{X_T}}$.
\end{cor}

\begin{proof}
By Lemma \ref{XtminusXs} we have that $X_T \in \bigcap_{p \in [1,\infty)}L_p$, so by Lemma \ref{Guniformboundforphi} we can choose 
$\varphi_\theta^{X_T,X_T^\pi} = \varphi_{C,\theta}^{X_T}$, where $\varphi_{C,\theta}^{X_T}$ is a bump function with decay faster than any polynomial. 
Now using Theorem \ref{Gpthm} 
we get for any $q \in [1,\infty)$ and $\theta \in (0,1)$ that
\begin{eqnarray*}
&&\E|g(X_T)-g(X_T^\pi)|^p \\
&&\leq  3 \cdot 2^p (\sup f_{X_T})^\frac{q(1-\theta)}{q+1} \left(V_{p,\varphi_{C,\theta}^{X_T}}(g)\right)^p C_q^\frac{q(1-\theta)}{q+1} \mesh^\frac{\gamma q(1-\theta)}{q+1}.
\end{eqnarray*}

Let $ 0 < \varepsilon < \gamma$. Now choose $q = \frac{2\gamma}{\varepsilon}-1$ and let $\theta=1/q$. Note that $q > 1$ since $\varepsilon < \gamma$. Then
$$ \frac{q(1-\theta)}{q+1}=\frac{q-1}{q+1} = 1 - \frac{\varepsilon}{\gamma} $$
and 
thus we get for all $g \in \G_{p,\varphi_{C,\theta}^{X_T}}$ that
$$ 
\E|g(X_T)-g(X_T^\pi)|^p 
\leq 3 \cdot 2^p (\sup f_{X_T})^{1 - \frac{\varepsilon}{\gamma}} \left(V_{p,\varphi_{C,\theta}^{X_T}}(g)\right)^p C_{1/\theta}^{1 - \frac{\varepsilon}{\gamma}} \mesh^{\gamma - \varepsilon}.
$$
Moreover, by Theorem \ref{polynomialsinG} we have that $\mathcal{P} \subset \G_{p,\varphi_{C,\theta}^{X_T}}$.
\end{proof}

\begin{remark}
In Corollary \ref{ExtensionRateCorGeneral} the function $\varphi_\theta^{X_T,X_T^\pi}$ depends on the distribution of $X_T^\pi$ and is replaced by the uniform bound $\varphi_{C,\theta}^{X_T}$. However, when considering convergence rate we are looking at partitions with small mesh size. Thus if approximating random variables $X_T^\pi$ corresponding to partitions with large mesh size had heavy tailed distributions, the use of the uniform bound could unnecessarily narrow down the class of functions. Therefore in such a case it would be better to take more delicate approach and study the result
$$ \E|g(X_T)-g(X_T^\pi)|^p \leq
  3 \cdot 2^p (\sup f_{X_T})^{1-\frac{\varepsilon}{\gamma}} \left(V_{p,\varphi_\theta^{X_T,X_T^\pi}}(g)\right)^p  
  C_{1/\theta}^{1-\frac{\varepsilon}{\gamma}} \mesh^{\gamma - \varepsilon}.$$
\end{remark}

Corollary \ref{ExtensionRateCorGeneral} now gives convergence rates for both Euler and Milstein schemes:

\begin{cor}
\label{ExtensionRateCorEulerMilstein}
Let $1 \leq p < \infty $. Then for $0 < \varepsilon < \half$, $\theta = \frac{\varepsilon}{1 - \varepsilon}$ and $g \in \G_{p,\varphi_{C,\theta}^{X_T}}$ that
$$ \E|g(X_T)-g(X_T^E)|^p \leq 3 \cdot 2^p (\sup f_{X_T})^{1 - 2\varepsilon} \left(V_{p,\varphi_{C,\theta}^{X_T}}(g)\right)^p C_{1/\theta}^{1 - 2\varepsilon} \mesh^{\half - \varepsilon}.$$
and similarly for $0 < \varepsilon < 1$ and $\theta = \frac{\varepsilon}{2 - \varepsilon}$  we have that
$$ \E|g(X_T)-g(X_T^M)|^p \leq 3 \cdot 2^p (\sup f_{X_T})^{1 - \varepsilon} \left(V_{p,\varphi_{C,\theta}^{X_T}}(g)\right)^p C_{1/\theta}^{1 - \varepsilon} \mesh^{1 - \varepsilon}.$$
Especially, the statements hold for any $g \in \Pcal$.
\end{cor}

\begin{example}
Let us generate a jump function by choosing the measure $\mu$ to be a sum of Dirac measures,
$$\mu = \sumZ{\frac{\alpha_k}{\varphi_{C,\theta}^{X_T}(a_k)} \delta_{a_k}},$$
where $\alpha_k,a_k \in \R$ for all $k\in \Z$ and $a_k \neq a_l$ for $k \neq l$. Then from the integrability condition for $\mu$ we see that $g \in \G_{p,\varphi_{C,\theta}^{X_T}}$ if
\begin{equation}
\label{jumps}
\sumZ{|\alpha_k| \varphi_{C,\theta}^{X_T}(a_k)}^\frac{1}{p} < \infty.
\end{equation}
Therefore the result of Corollary \ref{ExtensionRateCorGeneral} holds for jump functions with jumps controlled by the decay of the function $\varphi_{C,\theta}^{X_T}$ in a way that the condition \eqref{jumps} is satisfied.
\end{example}

\subsection{Euler scheme}

In the case of the Euler scheme we can again use our knowledge about constants to get more explicit results for the decay of the function $\varphi_\theta^{X_T,X_T^E}$. Let us recall the following result from Bouleau and L\'epingle \cite{BL}:

\begin{lemma}[{\cite[Ch. 5, Lemma B.1.2.]{BL}}]
\label{eulerpthmoment}
For $1\leq p < \infty$,
$$ \sup_{t \leq T} |X_t^E| \in L^p$$
and there exist $M(x_0,T, C_T) >0$ such that
$$ \norm{\sup_{t \leq T} |X_t^E|}{p} \leq e^{Mp^2}.$$
\end{lemma}

\begin{proof}
By Theorem \ref{Eulerthm} and Lemma \ref{XtminusXs} in the Appendix we get

$$
\begin{aligned}
 \norm{\sup_{t \leq T} |X_t^E|}{p} &\leq \norm{\sup_{t \leq T} |X_t - X_t^E|}{p} + \norm{\sup_{t \leq T} |X_t-x_0|}{p}+|x_0| \\
&\leq (2\sqrt{T}+|x_0|) e^{Mp^2},
\end{aligned}
$$
and we absorb the constant $(2\sqrt{T}+|x_0|)$ into the constant $M$.
\end{proof}

\begin{thm}
\label{varphiexpdecay}
We can choose the function $\varphi_\theta^{X_T,X_T^E}$ in a way that $\varphi_\theta^{X_T,X_T^E} \leq  \varphi_{\theta,E}^{X_T}$, where $\varphi_{\theta,E}^{X_T}$ is a bump function such that
\begin{enum_i}
\item if the functions $\sigma$ and $b$ are bounded, i.e. $|\sigma|,|b| \leq M$, we have
$$\varphi_{\theta,E}^{X_T}(z) = 
\begin{cases}
e^{-\frac{\theta(z-(x_0+MT))^2}{2M^2T}} & \text{ if } z > \max(x_0+MT,0),  \\
e^{-\frac{\theta(z-(x_0-MT))^2}{2M^2T}} & \text{ if } z < \min(x_0-MT,0),  \\
1 & \text{ elsewhere, } 
\end{cases}
$$
\item if the functions $\sigma$ and $b$ are Lipschitz, then there exists $z_0 > 1$ such that we have 
$$\varphi_{\theta,E}^{X_T}(z) =  
\begin{cases}
|z|^{-\frac{2\theta}{3\sqrt{3M}}(\log |z|)^{1/2}} & \text{ if } |z| > z_0, \\
1 & \text{ if } |z| \leq z_0, 
\end{cases}
$$
where $M=M(x_0,T,C_T)>0$.
\end{enum_i}
\end{thm}

\begin{proof}
\comment{
Since probabilities of the type $\Prob(X \leq -\lambda)$ are contained in probabilities $\Prob(|X| \geq \lambda)$, it suffices to consider only terms of the latter type.
Look first at the term $ \Prob(|X_T^\pi| \geq \lambda)$.
}

(i) We consider the Euler approximation with $n$ time nodes in the integral form \eqref{Eulercontinuousscheme}.
If we denote $$L_u:= \sum_{k=0}^{n-1} \sigma(t_k,X_{t_k}^E) \chi_{[t_k,t_{k+1})}(u),$$ then by the boundedness of $\sigma$ and the Novikov condition
$$M_t := e^{\alpha \int_0^t L_u \, dW_u - \frac{\alpha^2}{2} \int_0^t L_u^2 \, du} $$ is a martingale for any $\alpha > 0$, and $\E M_t = 1$.
Thus by Chebychev's inequality we have for $\lambda > 1$ that
$$ \Prob\left(e^{\alpha \int_0^T L_u \, dW_u - \frac{\alpha^2}{2} \int_0^T L_u^2 \, du} \geq \lambda\right) \leq \frac{1}{\lambda}. $$ By taking logarithm this implies 
$$ \Prob\left( \alpha \int_0^T L_u \, dW_u - \frac{\alpha^2}{2} \int_0^T L_u^2 \, du \geq \lambda\right) \leq e^{-\lambda} $$
for $\lambda > 0$. Since $$ \int_0^T L_u^2 \, du \leq M^2T,$$
we get 
$$ \Prob\left( \int_0^T L_u \, dW_u \geq \frac{\lambda}{\alpha} +\frac{\alpha M^2T}{2}\right) \leq e^{-\lambda}, $$
which we can reparametrize to get
$$ \Prob\left( \int_0^T L_u \, dW_u \geq \lambda \right) \leq e^{\frac{\alpha^2M^2T}{2} - \lambda \alpha} $$
for $\lambda > \alpha M^2T/2$.
Now we can choose $\alpha= \lambda/(M^2T)$ to get
$$ \Prob\left( \int_0^T L_u \, dW_u \geq \lambda \right) \leq e^{ -\frac{\lambda^2}{2M^2T}}$$
for $\lambda > 0$.
\comment{ auki laskettu versio
Now $\alpha$ is a free parameter, so we can choose $\alpha= c \lambda$ to get
$$ \Prob\left( \int_0^T L_u \, dW_u \geq \lambda \right) \leq e^{ \lambda^2\left( \frac{c^2M^2T}{2} - c\right)}. $$
The factor in the exponent is minimal if $c = 1/(M^2T)$, so for $\lambda > 0$ we have
$$ \Prob\left( \int_0^T L_u \, dW_u \geq \lambda \right) \leq 2e^{ -\frac{\lambda^2}{2M^2T}}.$$
}
A similar proof with $\widetilde L_u=-L_u$ shows that
$$ \Prob\left( \int_0^T L_u \, dW_u \leq \lambda \right) \leq e^{ -\frac{\lambda^2}{2M^2T}}$$
for $\lambda < 0$. Therefore, for $\lambda > \max(x_0+MT,0)$
\begin{eqnarray*}
\Prob\left(X_T^E \geq \lambda \right) &\leq& 
\Prob\left( x_0+\int_0^T L_u \, dW_u+MT \geq \lambda  \right) \\ 
&\leq& e^{ -\frac{(\lambda-(x_0+MT))^2}{2M^2T}},
\end{eqnarray*}
and for $\lambda < \min(x_0-MT,0)$
\begin{eqnarray*}
\Prob\left(X_T^E \leq \lambda \right) &\leq& 
\Prob\left( x_0+\int_0^T L_u \, dW_u-MT \leq \lambda  \right) \\ 
&\leq& e^{ -\frac{(\lambda-(x_0-MT))^2}{2M^2T}}.
\end{eqnarray*}
Obviously a similar proof works for the random variable $X_T$ instead of $X_T^E$, 
so by the definition of $\varphi_\theta^{X_T,X_T^E}$ the assertion follows. Moreover, to get a bump function we choose the upper bound to be one on the interval $[\min(x_0-MT,0),\max(x_0+MT,0)]$.

(ii)
If $\sigma$ and $b$ are Lipschitz, then we know from Lemma \ref{eulerpthmoment} that
$$\norm{X_T^E}{p} \leq e^{Mp^2},$$
where the constant $M>0$ depends on $x_0$, $T$ and $C_T$. 
Now by Chebychev's inequality we have for $\lambda > 0$ that
$$ \Prob(|X_T^E| \geq \lambda) \leq \frac{\E|X_T^E|^p}{\lambda^p} 
\leq \frac{e^{Mp^3}}{\lambda^p}.
$$
Choose $3Mp^2= \log \lambda$ for $\lambda > \lambda_0=e^{3M}$. This gives 
$$p=\frac{(\log \lambda)^{1/2}}{(3M)^{1/2}},$$ 
and thus for $\lambda > \lambda_0$ we get
$$ 
\Prob(|X_T^E| \geq \lambda) 
=  \frac{e^{\frac{1}{3}p\log \lambda}}{\lambda^p}=\lambda^{-\frac{2}{3}p} =
\lambda^{-\frac{2}{3\sqrt{3M}}(\log \lambda)^{1/2}}.
$$
Again the same proof works for the term $\Prob(|X_T| \geq \lambda)$ because of Lemma \ref{XtminusXs} in the Appendix.
\end{proof}

\begin{thm}
Let $c>0$. If the functions $\sigma$ and $b$ are bounded, and
$$g(z)=e^{c|z|^\gamma}\ \text{ if }\ 0< \gamma < 2,$$ 
or
$$g(z)=e^{c|z|^2} \text{ with } c < \theta/p,$$
then $g \in \G_{p,\varphi_{\theta,E}^{X_T}}$.
\end{thm}

\begin{proof}
Since $g$ is not differentiable at zero, define $\tilde g(x) := g'(x)$ if $x \neq 0$ and $\tilde g(0):=0$. By choosing a signed measure
$$d\mu(z)=\frac{\sgn(z) \tilde g(z)}{\varphi_{\theta,E}^{X_T}(z)}dz,$$
we get that $|\mu|$ is $\sigma$-finite and the representation $g=g(0)+g^\mu$ holds.
The definition of the class $\G_{p,\varphi_{\theta,E}^{X_T}}$ gives the condition
$$\intR \left(\varphi_{\theta,E}^{X_T}(z)\right)^\frac{1}{p} |\tilde g(z)| \, dz = \intR \left(\varphi_{\theta,E}^{X_T}(z)\right)^\frac{1}{p} e^{c|z|^\gamma} c \gamma |z|^{\gamma-1} \, dz < \infty,$$
which is by Theorem \ref{varphiexpdecay} satisfied, because the singularity at zero for $0 < \gamma < 1$ is not too strong, and  integrability is determined by the parameters $\gamma$, $c$, $\theta$ and $p$ as proposed in the formulation of this Theorem. Similarly we see that 
the local integrability condition is satisfied.
\end{proof}

\section{Lower bound}
\label{SectionLowerBound}

In this section we find a solution $X_1$ (i.e. $T=1$) of an SDE of the type \eqref{SDE} such that it gives a lower bound for the approximation rate of the Euler scheme in Theorem \ref{IndicatorRateThmEulerMilstein}. This is achieved by choosing $X_t=S_t$, the geometric Brownian motion. Let $S_t=e^{W_t-t/2}$ for $t \in [0,1]$, so that $S$ is a solution of 
$$ S_t=1+\int_0^t S_s \, dW_s $$
and let $U^n:= S^E - S$, where $S^E$ is the Euler scheme as defined in (\ref{Eulercontinuousscheme}) corresponding to the equidistant partition of $[0,1]$, i.e. $\pi=(i/n)_{i=0}^n$. 

\begin{lemma}
\label{JPlemma}
We have $(W,\sqrt{n}U^n) \Longrightarrow (W,U)$ in the Skorohod topology, where $U$ is the strong solution of the equation
\begin{equation}
\label{JPlemmaeq}
 U_t=\int_0^t U_s \, dW_s - \frac{1}{\sqrt{2}}\int_0^t S_s \, dB_s
\end{equation}
and $B$ is a standard Brownian motion independent of $W$. 
\end{lemma}
\begin{proof}
The statement is an immediate consequence of a result by Jacod and Protter, \cite[Corollary 5.4]{JP}.
\end{proof}

\begin{thm}
\label{lowerboundthm}
There exists $K_0 > 0$ such that
$$\liminf_{n \to \infty} \sqrt{n} \sup_{K \geq K_0} \E|\chi_{[K, \infty)}(S_1)-\chi_{[K, \infty)}(S_1^E)| > 0, $$
where $S_1^E$ is the equidistant Euler approximation of $S_1$.
\end{thm}

\begin{remark}
Theorem \ref{lowerboundthm} states that the convergence rate $\half-\varepsilon$ for the Euler scheme obtained in Theorem \ref{IndicatorRateThmEulerMilstein} and consequently in Theorem \ref{BVRateThmEulerMilstein} and Corollary \ref{ExtensionRateCorEulerMilstein} is optimal up to the factor $\varepsilon$, i.e. any rate $\gamma > \half$ leads to a contradiction with the statement of Theorem \ref{lowerboundthm}.
\end{remark}

\begin{proof}[Proof of Theorem \ref{lowerboundthm}]
\comment{ Since $W$ is a continuous local martingale with $C_t=t$, choose $Y=W$ in Lemma \ref{JPlemma}. Then let $\sigma(x)=x$, which is in $C^1$ with at most linear growth, and $x_0=1$. Thus the equation is
$$ X_t=1 + \int_0^t X_s\ dW_s. $$
The solution is the geometric Brownian motion, which has a bounded density. Then}
Let us consider the setting of Lemma \ref{JPlemma} and the process $U$ defined by the equation \eqref{JPlemmaeq}. If $U_1=0$ a.s., then for all $t \in [0,1]$ we have $U_t=0$ a.s., which leads to a contradiction. 
Therefore $\Prob\left(U_1 > 0\right) > 0$ or $\Prob\left(U_1 < 0\right) > 0$. If $\Prob\left(U_1 > 0\right) > 0$, then there exist $\varepsilon \in (0,1]$, $\delta > 0$ and $K \geq 1+K_0$ with $K_0 > 0$ such that
$$\Prob\left( S_1 \in [K-1,K),\ U_1 > \varepsilon \right) = \delta.$$
The case $\Prob\left(U_1 < 0\right) > 0$ can be treated in a similar way by changing the condition $U_1 > \varepsilon$ to $U_1 < -\varepsilon$.
By Lemma \ref{JPlemma} we know that $(W,\sqrt{n}U^n) \Rightarrow (W,U)$ in the Skorohod topology. This implies that  $(W_1,\sqrt{n}U^n_1) \Rightarrow (W_1,U_1)$, since the projection mapping $\pi_1$, i.e. the mapping $\alpha \mapsto \alpha(1)$ for a process $\alpha$, is continuous in the Skorohod topology. 
Because the function $e^{x-\frac{t}{2}}$ is continuous, we have $(S_1,\sqrt{n}U^n_1) \Rightarrow (S_1,U_1)$. Therefore
$$
\begin{aligned}
&\liminf_{n \to \infty} \Prob\left( S_1 \in [K-1,K), \sqrt{n}[S_1^E - S_1] > \varepsilon \right) \\
&= \liminf_{n \to \infty} \Prob\left( S_1 \in (K-1,K), \sqrt{n}\,U_1^n  > \varepsilon \right) \\
&\geq\Prob\left( S_1 \in (K-1,K), U_1 > \varepsilon \right) \\
&=\Prob\left( S_1 \in [K-1,K), U_1 > \varepsilon \right),
\end{aligned}
$$
and we see that there exists $n_0 \geq 1$ such that for all $n \geq n_0$
$$ \Prob\left( S_1 \in [K-1,K),\ [S_1^E - S_1] > \frac{\varepsilon}{\sqrt{n}} \right) \geq \frac{\delta}{2}.$$
Assume a partition $K-1=K_0^m < K_1^m < \cdots < K_m^m = K$. Then
$$ \sup_{l=1,\dots,m} \Prob\left( S_1 \in [K_{l-1}^m,K_l^m),\ [S_1^E - S_1] > \frac{\varepsilon}{\sqrt{n}} \right) \geq \frac{\delta}{2m}.$$
Now choose the partition $(K_l^m)_{l=1}^m$ to be equidistant with
\begin{equation}
\label{meshsizecondition}
 \frac{1}{m} \leq \frac{\varepsilon}{\sqrt{n}}. 
\end{equation}
Then there exists $l_0 \in \{1, \dots , m\}$ such that 
$$
\begin{aligned}
\frac{\delta}{2m} &\leq \Prob\left( S_1 \in [K_{l_0-1}^m,K_{l_0}^m),\  S_1^E  >  S_1 + \frac{\varepsilon}{\sqrt{n}} \right) \\
&\leq \Prob\left( S_1 < K_{l_0}^m,\  S_1^E  \geq  K_{l_0}^m \right).
\end{aligned}
$$
Let $m = \lceil\sqrt{n}/\varepsilon \rceil$, which satisfies the condition \eqref{meshsizecondition} for the mesh size. Hence
\begin{eqnarray*}
\frac{\delta}{2\lceil\sqrt{n}/\varepsilon \rceil} 
&\leq& \Prob\left( S_1 < K_{l_0}^m,\  S_1^E  \geq  K_{l_0}^m \right) \\
&\leq& \E|\chi_{[K_{l_0}^m, \infty)}(S_1)-\chi_{[K_{l_0}^m, \infty)}(S_1^E)|.
\end{eqnarray*}
Since $\lceil \sqrt{n}/ \varepsilon\rceil \leq 2\sqrt{n} / \varepsilon $ we have
$$
\E|\chi_{[K_{l_0}^m, \infty)}(S_1)-\chi_{[K_{l_0}^m, \infty)}(S_1^E)|
\geq \frac{\delta}{2\lceil\sqrt{n}/ \varepsilon \rceil} \geq \frac{\delta \varepsilon}{4\sqrt{n}}.  
$$
Therefore
$$
\sqrt{n}\sup_{K \geq K_0} \E|\chi_{[K, \infty)}(S_1)-\chi_{[K, \infty)}(S_1^E)|
\geq \frac{\delta \varepsilon}{4}  
$$
for all $n \geq n_0$, which implies the assertion. 
\comment{
Moreover, it is easy to see that 
$$S^E_1=\prod_{i=1}^n (1+W_{t_i}-W_{t_{i-1}}), $$
which is a product of independent Gaussian random variables. This implies that there exists a constant $\delta_{n_0} > 0$ such that 
$$\inf_{1 \leq n < n_0} \Prob (S_1^E < 0) \geq \delta_{n_0}.$$
Thus we have
\begin{eqnarray*}
&&\!\!\!\sqrt{n}\sup_{K \geq 0} \E|\chi_{[K, \infty)}(S_1)-\chi_{[K, \infty)}(S_1^E)| 
\geq \sqrt{n}\E|\chi_{[0, \infty)}(S_1)-\chi_{[0, \infty)}(S_1^E)| \\
&&\geq \sqrt{n}\E|1-\chi_{[0, \infty)}(S_1^E)| 
= \sqrt{n}\Prob (S_1^E < 0) \geq \delta_{n_0} 
\end{eqnarray*}
for all $1 \leq n < n_0$.
Thus the assertion is true for all $n \geq 1$.
}
\end{proof}

\appendix

\section{}
\label{SectionAppendix}

Here we prove the following Theorem from the book of Bouleau and L\'epingle, \cite[pp. 275-276]{BL}. The proof is given in the book, but without computing the constant explicitly.

\begin{thm} 
\label{Eulerthm}
If the assumptions (i)-(iii) in section \ref{SectionAssumptions} hold, and $1 \leq p < \infty$, then
$$  \norm{\sup_{0 \leq t \leq T} |X_t-X_t^E|}{p} \leq e^{Mp^2} \mesh^\half,$$
where the constant $M > 0$ depends at most on $x_0$, $T$ and $C_T$.
\end{thm}

For the proof we need the following Lemma:

\begin{lemma}
\label{XtminusXs}
For $1\leq p < \infty$ and $0 \leq s \leq t\leq T$,
$$\norm{\sup_{u \in [s,t]} |X_u-X_s|}{p} \leq \sqrt{t-s}\,e^{Mp^2},$$
where $M > 0$ depends at most on $x_0$, $T$ and $C_T$.
\end{lemma}
\begin{proof}
Without loss of generality we can suppose that $p \geq 2$. For fixed $s \in [0,T]$ define 
$$F(t)=\norm{\sup_{u \in [s,t]}|X_u-X_s|}{p} $$
for all $t \in [s,T]$. Then by the Burkholder-Davis-Gundy inequality (as in \cite[p. 269]{BL}) 
and the linear growth condition we get
\begin{eqnarray*}
F(t)
&\leq& 2\norm{\int_s^t \sigma(u,X_u)\ dW_u}{p}+\norm{\int_s^t |b(u,X_u)|\, du}{p} \\
&\leq& 8p\norm{\left(\int_s^t \! |\sigma(u,X_u)|^2 \, du\right)^\half}{p}\!\!\! +\sqrt{t-s}\norm{\! \left(\int_s^t\! |b(u,X_u)|^2\, du\right)^\half}{p} \\
&\leq& C_T(8p+\sqrt{t-s})\norm{\left(\int_s^t (1+|X_u|)^2 \, du\right)^\half}{p} \\
&\leq& \widetilde{C}\left[\norm{\left(\int_s^t (1+|X_s|)^2 \, du\right)^\half}{p}+\norm{\left(\int_s^t (|X_u-X_s|)^2 \, du\right)^\half}{p}\right] \\
&\leq& \widetilde{C}\left[\sqrt{t-s}\norm{1+|X_s|}{p}  +\norm{\left(\int_s^t \sup_{v \in [s,u] }|X_v-X_s|^2 \, du\right)^\half}{p}\right] \\
&\leq& \widetilde{C}\left[\sqrt{t-s}(1+\norm{X_s}{p})  +\left(\int_s^t \norm{\sup_{v \in [s,u] }|X_v-X_s|}{p}^2 \, du\right)^\half\right] \\
&=& \widetilde{C}\sqrt{t-s}(1+\norm{X_s}{p}) +\widetilde{C}\left(\int_s^t F(u)^2 \, du\right)^\half,
\end{eqnarray*}
where $\widetilde{C}=C_T(8p+\sqrt{t-s})$. Thus we have
\begin{eqnarray*}
F(t)^2 &\leq& 2C_T^2(8p+\sqrt{t-s})^2(t-s)(1+\norm{X_s}{p})^2 \\
& & +2C_T^2(8p+\sqrt{t-s})^2\int_s^t F(u)^2 \, du
\end{eqnarray*}
and by the Gronwall lemma we get
$$F(t)^2 \leq 2C_T^2(8p+\sqrt{t-s})^2(t-s)(1+\norm{X_s}{p})^2e^{2C_T^2(8p+\sqrt{t-s})^2(t-s)}. $$
By taking the square root and choosing a suitable constant $M$ we get
$$F(t) \leq Mp\sqrt{t-s}(1+\norm{X_s}{p})e^{Mp^2}, $$
where $M=M(T,C_T)$. In particular, the above estimate gives that
$$\norm{X_s}{p} \leq |x_0| +\norm{X_s - X_0}{p}  \leq |x_0|+ Mp\sqrt{s}(1+|x_0|)e^{Mp^2}, $$
so if we 
redefine the constant $M$, we get
$$\norm{\sup_{u \in [s,t]}|X_u-X_s|}{p} \leq \sqrt{t-s} \,e^{Mp^2}, $$
where $M=M(x_0,T,C_T)$.
\comment{ tämä tuli käsiteltyä alussa (wlog assume p \geq 2)
The case $1 \leq p < 2$ follows from the case $p=2$, since 
$$\norm{\sup_{u \in [s,t]}|X_u-X_s|}{p} \leq \norm{\sup_{u \in [s,t]}|X_u-X_s|}{2} \leq \sqrt{t-s} \,e^{4M}.  $$
}
\end{proof}

\begin{proof}[Proof of Theorem \ref{Eulerthm}]

Suppose that $p \geq 2$. Now define 
$$ F(t) = \norm{\sup_{s \leq t} |X_s-X_s^E|}{p}. $$
Here $X^E$ is the Euler scheme related to the equidistant partition $\pi$, and is defined for continuous time by formula \eqref{Eulercontinuousscheme}.
Then by the Burkholder-Davis-Gundy inequality (\cite[p. 269]{BL}) we get
\begin{eqnarray*}
F(t) 
&\leq& 8p \norm{ \left( \int_0^t \sum_{k=0}^{n-1} |\sigma(u,X_u)- \sigma(t_k,X_{t_k}^E)|^2 \chi_{[t_k,t_{k+1})}(u)\, du \right)^\half}{p} \\
&& + \sqrt{t} \norm{ \left(\int_0^t \sum_{k=0}^{n-1} |b(u,X_u)- b(t_k,X_{t_k}^E)|^2 \chi_{[t_k,t_{k+1})}(u)\, du \right)^\half}{p} \\
&\leq& 8p \left( \int_0^t \left(\sum_{k=0}^{n-1} \norm{ \sigma(u,X_u)- \sigma(t_k,X_{t_k}^E)}{p}\chi_{[t_k,t_{k+1})}(u) \right)^2 \, du \right)^\half \\
&& + \sqrt{t} \left(\int_0^t \! \left( \sum_{k=0}^{n-1} \norm{ b(u,X_u)- b(t_k,X_{t_k}^E) }{p}\chi_{[t_k,t_{k+1})}(u)\right)^2\, du \right)^\half\!\!. 
\end{eqnarray*}
By the 
conditions in Section \ref{SectionAssumptions} we have that
\begin{eqnarray*}
&&|\sigma(u,X_u) - \sigma(t_k,X_{t_k}^E)| \leq |\sigma(u,X_u)- \sigma(u,X_{t_k})| \\
&& \ + |\sigma(u,X_{t_k})- \sigma(t_k,X_{t_k})|
 + |\sigma(t_k,X_{t_k})- \sigma(t_k,X_{t_k}^E)| \\
&&\leq C_T\left( |X_u-X_{t_k}| + (1+|X_{t_k}|)|u-t_k|^\alpha + |X_{t_k}-X_{t_k}^E|  \right), 
\end{eqnarray*}
and a corresponding inequality holds for the function $b$. Thus by Lemma \ref{XtminusXs} 
we have for $u \geq t_k$ that
\begin{eqnarray*}
&&\norm{\sigma(u,X_u) - \sigma(t_k,X_{t_k}^E)}{p} \\
&&\leq C_T \left( \norm{X_u-X_{t_k}}{p}+ \norm{(1+|X_{t_k}|)|u-t_k|^{\alpha }}{p} + \norm{X_{t_k}-X_{t_k}^E}{p} \right) \\
&&\leq C_T\left( e^{Mp^2}|u-t_k|^\half + \left(1+|x_0|+\sqrt{t_k}e^{Mp^2}\right)|u-t_k|^\alpha + F(u)  \right) \\
&&\leq C_T\left((1+|x_0|+\sqrt{T})e^{Mp^2}(|\pi|^\half + |\pi|^\alpha) + F(u)\right),
\end{eqnarray*}
and again a corresponding inequality holds for $b$. Denote $C(x_0,T)=1+|x_0|+\sqrt{T}$. Now we can continue our estimate for $F(t)$ to get
\begin{eqnarray*}
&&\!\!\! F(t) \leq C_T(8p+\sqrt{T})\ \cdot \\
&&\!\!\! \left( \int_0^t \left[\sum_{k=0}^{n-1}\left[ C(x_0,T)e^{Mp^2}(|\pi|^\half + |\pi|^\alpha) + F(u)\right]\chi_{[t_k,t_{k+1})}(u) \right]^2 \, du \right)^\half \\
&& \!\!\! \leq C_T(8p+\sqrt{T}) \left( \int_0^t \left(C(x_0,T)e^{Mp^2}(|\pi|^\half + |\pi|^\alpha) + F(u) \right)^2 \, du \right)^\half \\
&&\!\!\! \leq  C_T(8p+ \sqrt{T}) \left[ \sqrt{T} \,C(x_0,T)e^{Mp^2}(|\pi|^\half\! + \! |\pi|^\alpha)
+\left[\int_0^t F(u)^2\, du \right]^\half \right] \\
&& \!\!\! \leq e^{M_1 p^2}(|\pi|^\half + |\pi|^\alpha) + M_1 p \left(\int_0^t F(u)^2\, du \right)^\half,
\end{eqnarray*}
where $M_1=M_1(x_0,T,C_T)$ does not depend on $p$. Taking a square we get
$$
F(t)^2 \leq 2e^{2M_1 p^2}(|\pi|^\half + |\pi|^\alpha)^2 + 2M_1^2 p^2 \int_0^t F(u)^2\, du,
$$
and thus Gronwall's Lemma gives
$$F(t)^2 \leq 2e^{2M_1 p^2}(|\pi|^\half + |\pi|^\alpha)^2 e^{2M_1^2 p^2 T}=2e^{(2M_1+2M_1^2T) p^2}(|\pi|^\half + |\pi|^\alpha)^2. $$
By taking square root on both sides and recalling the assumption $\alpha \geq \half$ this gives
$$F(t) \leq \sqrt{2}e^{(M_1+M_1^2T) p^2}(|\pi|^\half + |\pi|^\alpha) \leq e^{M_2 p^2}|\pi|^\half, $$
where $M_2=M_2(x_0,T,C_T)$.
The case $1 \leq p < 2$ follows from the case $p=2$ by redefining the constant $M_2$.
\comment{ vanha pitempi loppu
By taking square root on both sides this gives
$$F(t) \leq \sqrt{2}e^{(M_1+M_1^2T) p^2}(|\pi|^\half + |\pi|^\alpha). $$
Set $\beta = \alpha \wedge \half$ and suppose that $|\pi|$ is small. Then 
$$F(t)= \norm{\sup_{s \leq t} |X_s-X_s^E|}{p} \leq e^{M_2 p^2}|\pi|^\beta, $$
where $M_2=M_2(x_0,T,C_T)$. We assume that $\alpha \geq \half$, so the convergence rate is $\beta=\half$. 
The case $1 \leq p < 2$ follows from the case $p=2$ by redefining the constant $M_2$.
}
\end{proof}


\begin{thebibliography}{0}
{ 

\bibitem{BL}
\textsc{Nicolas Bouleau, Dominique L\'epingle},
\textit{Numerical Methods for Stochastic Processes}.
Wiley, 1994.

\bibitem{BT1}
\textsc{Vlad Bally, Denis Talay},
\textit{The Law of the Euler Scheme for Stochastic Differential Equations: I. Convergence Rate of the Distribution Function}.
Probab. Theory Related Fields 104 (1996), no. 1, 43--60.

\bibitem{BT2}
\textsc{Vlad Bally, Denis Talay},
\textit{The Law of the Euler Scheme for Stochastic Differential Equations: II. Convergence Rate of the Density}.
Monte Carlo Methods Appl. 2 (1996), no. 2, 93--128.

\bibitem{BS}
\textsc{Colin Bennett, Robert Sharpley},
\textit{Interpolation of Operators}.
Academic Press, 1988.

\bibitem{CFN}
\textsc{María Emilia Caballero, Begoña Fernández, David Nualart},
\textit{Estimation of Densities and Applications}.
J. Theoret. Probab. 11 (1998), no. 3, 831 -- 851.

\bibitem{F}
\textsc{Avner Friedman},
\textit{Partial Differential Equations of Parabolic Type}.
Prentice-Hall, 1964.

\bibitem{HMR1}
\textsc{Norbert Hofmann, Thomas Müller-Gronbach, Klaus Ritter},
\textit{The Optimal Discretization of Stochastic Differential Equations}.
J. Complexity 17 (2001), No.1, 117 -- 153.

\bibitem{HMR2}
\textsc{Norbert Hofmann, Thomas Müller-Gronbach, Klaus Ritter},
\textit{Linear vs. Standard Information for Scalar Stochastic Differential Equations}.
J. Complexity 18 (2002), 394 -- 414.

\bibitem{HM}
\textsc{Norbert Hofmann, Thomas Müller-Gronbach},
\textit{On the Global Error of Itô-Taylor Schemes for Strong Approximation of Scalar Stochastic Differential Equations}.
J. Complexity 20 (2004), 732 -- 752.

\bibitem{JP}
\textsc{Jean Jacod, Philip Protter},
\textit{Asymptotic Error Distributions for the Euler Method for Stochastic Differential Equations}.
Ann. Prob. 26 (1998),  no. 1, 267--307.

\bibitem{KS}
\textsc{Ioannis Karatzas, Steven E. Shreve},
\textit{Brownian Motion and Stochastic Calculus, Second Edition},
Springer-Verlag, 1991.

\bibitem{KP}
\textsc{Peter E. Kloeden, Eckhard Platen},
\textit{Numerical Solutions of Stochastic Differential Equations}.
Springer-Verlag, 1992.

\bibitem{M1}
\textsc{Thomas Müller-Gronbach},
\textit{Strong Approximation of Systems of Stochastic Differential Equations}.
Habilitation thesis. Darmstadt, 2002.

\bibitem{M2}
\textsc{Thomas Müller-Gronbach},
\textit{The Optimal Uniform Approximation of Systems of Stochastic Differential Equations}.
Ann. Appl. Probab. 12 (2002),  no.2, 664 -- 690.

\bibitem{M3}
\textsc{Thomas Müller-Gronbach},
\textit{Optimal Pointwise Approximation of SDEs Based on Brownian Motion at Discrete Points}.
Ann. Appl. Probab. 14 (2004),  no. 4, 1605--1642.

\bibitem{MSTZ}
\textsc{Kyoung-Sook Moon, Anders Szepessy, Raúl Tempone, Georgios E. Zouraris},
\textit{Convergence Rates for Adaptive Weak Approximation of Stochastic Differential Equations}.
Stoch. Anal. Appl.  23  (2005),  no. 3, 511--558.

\bibitem{N}
\textsc{David Nualart},
\textit{The Malliavin Calculus and Related Topics}.
Springer-Verlag, 1995.

\bibitem{WR}
\textsc{Walter Rudin},
\textit{Real and Complex Analysis, Second Edition}.
McGraw-Hill, 1966, 1974.



}

\end{thebibliography}
\end{document}